\documentclass[12pt,leqno]{article}
\title{\ }
\author{\ }
\begin{document}
\newtheorem{theorem}{Theorem}[section]
\newtheorem{proposition}{Proposition}[section]
\newtheorem{lemma}{Lemma}[section]
\newtheorem{definition}{Definition}[section]
\newtheorem{corollaire}{Corollary}[section]
\newtheorem{remarque}{Remark}[section]
\def\Lim{\displaystyle\lim}
\def\Sup{\displaystyle\sup}
\def\Inf{\displaystyle\inf}
\def\Max{\displaystyle\max}
\def\Min{\displaystyle\min}
\def\Sum{\displaystyle\sum}
\def\Frac{\displaystyle\frac}
\def\Int{\displaystyle\int}
\def\n{|\kern -.05cm{|}\kern -.05cm{|}}
\def\Bigcap{\displaystyle\bigcap}
\def\E{{\cal E}}
\def\R{{\bf \hbox{\sc I\hskip -2pt R}}} 
\def\N{{\bf \hbox{\sc I\hskip -2pt N}}} 
\def\Z{{\bf Z}} 
\def\Q{{\bf \hbox{\sc I\hskip -7pt Q}}} 
\def\C{{\bf \hbox{\sc I\hskip -7pt C}}} 

{\begin{center} {\ } \vskip 2.1cm {\large\bf
Indirect boundary stabilization of weakly coupled degenerate wave equations with a single fractional
damping at the degenerate boundary
 \\ } \
\\ \ \\
\begin{tabular}{c}
\sc Chahira Aichi and Abbes Benaissa,
\\ \\
\small  Laboratory of Analysis and Control of PDEs, \qquad \\
\small Djillali Liabes University,\qquad \\
\small P. O. Box 89, Sidi Bel Abbes 22000, ALGERIA.\qquad \\
\small E-mails: aichi.math@gmail.com  \\
\small benaissa$_{-}$abbes@yahoo.com \\
\end{tabular}
\end{center}
\
\\
\begin{abstract}
In this paper, we investigate the well-posedness and stability of a one-dimensional system of
degenerate wave equations coupled through zero-order terms and subject to a single boundary fractional damping mechanism
acting at the degenerate endpoint. We establish an optimal polynomial energy decay rate of order
$1/t^{(2+2\nu_{\gamma}-\tau)}$.
Our analysis is based on the frequency-domain approach and exploits an explicit Bessel-function representation of the resolvent operator.
\end{abstract}

{\it Keywords :} System of coupled degenerate wave equations, Fractional boundary damping, Strong
asymptotic stability, Bessel functions, Optimal polynomial decay.

{\it ${\cal AMS}$ Classification:} 35B40,35L80, 74D05, 93D15.
\begin{sloppypar}
\section{Introduction}
Coupled wave equations arise naturally in the mathematical
description of numerous physical phenomena involving the interaction
of elastic, mechanical, and acoustic structures (see {\bf\cite{lima}}, {\bf\cite{schm}} and {\bf\cite{zazaza}}).
Examples include coupled strings and beams, elastic networks, wave propagation in
composite media, and mechanical systems consisting of interconnected
components. In such models, the interaction between different
components is generally represented by coupling terms describing the
exchange of mechanical energy between the subsystems. Understanding
the long-time behavior of these systems, particularly their
stability and energy dissipation properties, constitutes an
important problem in the theory of evolution equations.

An additional difficulty arises when the propagation coefficient is
space-dependent and degenerates at a boundary point or within the
spatial domain. Degenerate hyperbolic equations have attracted
considerable attention because the degeneracy modifies the
propagation properties of waves and introduces analytical
difficulties that are not present in uniformly hyperbolic problems.
In particular, when the coefficient of the principal spatial
operator vanishes at an endpoint, the associated energy space is
naturally described by weighted Sobolev spaces. The degeneracy may
also affect the regularity of solutions, the behavior of boundary
traces, and the spectral properties of the underlying differential
operator.

In this work, we investigate the existence and energy decay rate of a system of
coupled degenerate wave equations, that is,
\begin{equation}
\left\{
\matrix{
u_{tt}(x,t)-(x^{\gamma}u_{x})_{x}(x,t)+\alpha v=0 \hfill & \hbox{ in }(0,1)\times
(0,+\infty ), \hfill & \cr
v_{tt}(x,t)-(x^{\gamma}v_{x})_{x}(x,t)+\alpha u=0\hfill & \hbox{ in }(0,1)\times
(0,+\infty ),\hfill & \cr
(x^{\gamma} u_{x})(0,t)=\rho \partial_{t}^{\tau, \omega} u(0, t)\hfill & \hbox{ in } (0, +\infty),\hfill&\cr
u(1,t)=v(0, t)=v(1,t)=0 \hfill & \hbox{ for } t\in (0, +\infty),\hfill & \cr
u(x,0)=u_0(x), u_{t}(x,0)=u_1(x)\hfill & \hbox{ for } x\in (0, 1),\hfill &\cr
v(x,0)=v_0(x), v_{t}(x,0)=v_1(x)\hfill & \hbox{ for } x\in (0, 1),\hfill &\cr}
\right.
\label{eqq100}
\end{equation}
where $\gamma\in [0, 1)$ and $\rho>0$. The notation
$\partial_{t}^{\tau, \omega}$ stands for the generalized Caputo's
fractional derivative of order $\tau$, $(0< \tau\leq 1)$, with
respect to the time variable (see {\bf\cite{choi}}). It is defined as follows
$$
\partial^{\tau,\omega}_{t}g(t)=\left\{
\matrix{
g_{t} & \hbox{for} \ \tau=1, \ \omega\geq 0,\hfill &\cr
\frac{1}{\Gamma(1-\tau)}\int_{0}^{t}(t-s)^{-\tau}e^{-\omega(t-s)}  \frac{dg}{ds}(s)ds\quad &
\hbox{for }\ 0<\tau <1,\quad\omega\geq 0. \ \hfill &\cr}
\right.
$$
The initial data $(u_{0},u_{1},v_{0},v_{1})$ belong to a suitable
function space.

The coefficient \(x^\gamma\) degenerates at the endpoint $x=0$, while the parameter $\alpha$
describes the interaction between the two wave components. The
system can be viewed as a pair of coupled elastic media whose
effective tension or stiffness decreases near one endpoint. The
coupling terms express the mechanical interaction between the two
components and allow energy to be exchanged between the
corresponding subsystems.

We begin by reviewing the existing literature on linear degenerate wave equations, with particular emphasis on both qualitative and quantitative aspects
of their analysis. In the scalar case, we first recall the pioneering work of
Alabau-Boussouira et al. \cite{alabau1}, which considers the problem on the domain $(x,t)\in(0,1)\times(0,+\infty)$, focusing primarily on the following equation:
$$
u_{tt}(x,t)-(a(x) u_{x}(x,t))_x=0 \hbox{ in } (0, 1)\times (0, \infty).
$$
The system is supplemented with the following boundary conditions, incorporating linear boundary damping:
$$
\left\{\matrix{\left\{\matrix{u(0, t)=0 \hfill & \hbox{ if } 0\leq m_a< 1 \hfill & \cr
(a u_{x})(0,t)=0 \hfill & \hbox{ if } 1\leq m_a< 2 \hfill & \cr}\right.\hfill  & \hbox{ in } (0, +\infty) ,\hfill&\cr
u_t(1, t)+u_x(1, t)+\beta u(1, t)=0 \hfill & \hbox{ in } (0, +\infty).\hfill&\cr
}\right.
$$
where $\beta>0$ is the given constant. $m_a=\sup_{0< x\leq 1}\Frac{x|a'(x)|}{a(x)}<2$
quantifies the degree of degeneracy of the diffusion coefficient.
By means of the energy multiplier method, the authors establish exponential decay of the total energy of the system.

Recently, in {\bf\cite{beai}}, Benaissa and Aichi studied the stabilization of a scalar degenerate wave equation
equipped with the following fractional boundary damping
$$
\left\{\matrix{\left\{\matrix{u(0, t)=0 \hfill & \hbox{ if } 0\leq m_a< 1 \hfill & \cr
(a u_{x})(0,t)=0 \hfill & \hbox{ if } 1\leq m_a< 2 \hfill & \cr}\right.\hfill  & \hbox{ in } (0, +\infty) ,\hfill&\cr
(au_{x})(1, t)+\rho\partial_{t}^{\tau, \omega}u(1, t)+\beta u(1,t)=0 \hfill & \hbox{ in } (0, +\infty).\hfill&\cr
}\right.
$$
They established optimal polynomial stability of the solutions by employing a frequency-domain
approach combined with suitable multiplier techniques.

In {\bf\cite{zerk}}, Zerkouk, Aichi, and Benaissa studied the stability of a one-dimensional weakly degenerate
wave equation subject to a fractional boundary feedback acting at the degenerate endpoint. By means of semigroup
theory and resolvent estimates, they established several stability properties, including strong and uniform stability,
as well as explicit polynomial decay estimates under suitable assumptions. Their work provides a useful framework for investigating
the influence of the interaction between the fractional boundary damping and degenerate coefficient on the long-time dynamics of degenerate wave equations.

Subsequently, Liu and Rao {\bf\cite{lira}} investigated a general class of coupled systems of second-order evolution equations.
The system is governed by the following equations:
$$
\left\{\matrix{u_{tt}-b\Delta u+\alpha y=0 \hfill & \hbox{ on } \Omega,\hfill \cr
y_{tt}-\Delta u+\alpha u=0 \hfill & \hbox{ on } \Omega,\hfill \cr
u=0 \hfill  &\hbox{ on } \Gamma_D, \hfill \cr
b\partial_{\nu} u+\gamma u+u_t=0\hfill &\hbox{ on }\Gamma_N,\hfill\cr
y=0 \hfill&\hbox{ in } \Gamma,\hfill\cr}\right.
$$
where $\Omega\subset \R ^n$ is a bounded domain with smooth boundary $\Gamma$ of class $C^2$
such that $\Gamma=\Gamma_D\cup \Gamma_N$ and $\Gamma_D\cap \Gamma_N=\emptyset$.

Using the frequency-domain approach, they established a polynomial energy decay rate of order $\frac{\ln t}{t}$ for sufficiently smooth initial data in
the equal wave-speed case. Furthermore, when the wave speeds
are different (i.e., $b\not= 1$), they showed that the energy decay rate is governed by the arithmetic properties of the ratio of the wave speeds.

More recently, Wehbe and Koumaiha {\bf\cite{weko}} investigated a one-dimensional system of wave equations coupled through zero-order terms. More precisely,
they addressed the stabilization problem for the following system of partially damped coupled wave equations propagating at equal speeds:
$$
\left\{
\matrix{
u_{tt}-u_{xx}+\alpha y=0 \hfill & \hbox{ in }(0,1)\times (0,+\infty ), \hfill & \cr
y_{tt}-y_{xx}+\alpha u=0\hfill & \hbox{ in }(0,1)\times (0,+\infty ),\hfill & \cr
u(0, t)=y(0, t)=y(1, t)=0 \hfill  & \hbox{ in } (0, +\infty) ,\hfill&\cr
u_{x}(1, t)+\rho u_t(1, t)=0 \hfill & \hbox{ in } (0, +\infty),\hfill&\cr
u(x,0)=u_0(x), u_{t}(x,0)=u_1(x),
y(x,0)=y_0(x), y_{t}(x,0)=y_1(x)\hfill & \hbox{ for } x\in (0, 1),\hfill &\cr}
\right.
$$
where $\rho> 0$. They established an optimal polynomial energy decay rate of order $\frac{1}{t}$
by employing the frequency-domain approach together with the Riesz basis property of the system's generalized eigenvectors.

In {\bf\cite{akil}}, Akil et al. investigated the indirect boundary stabilization of a one-dimensional system of coupled wave equations,
governed by the following model:
$$
\left\{
\matrix{
u_{tt}(x,t)-u_{xx} (x,t)-d v_t(x,t)=0 \hfill &\hbox{ in }
(0, 1)\times (0,+\infty ), \hfill & \cr
v_{tt}(x,t)-v_{xx}+d u_t(x,t)=0\hfill & \hbox{ in }%
(0, 1)\times (0,+\infty ),\hfill & \cr
u(0,t)=v(0,t)=v(1, t)=0 \hfill & \hbox{ on } (0,+\infty ),\hfill & \cr
u_{x}(1,t)+\rho \partial _{t}^{\tau ,\omega }u(1,t)=0
\hfill & \hbox{ on } (0,+\infty ),\hfill & \cr
u(x,0)=u_{0}(x),\ u_{t}(x,0)=u_{1}(x)\hfill & \hbox{ on }
(0,1),\hfill & \cr
v(x,0)=v_{0}(x),\ v_{t}(x,0)=v_{1}(x)\hfill & \hbox{ on }(0,1).\hfill & \cr}
\right.
$$
They proved that the energy decays polynomially at the rate $t^{-s(\tau)}$, with

\noindent
i) If $d\not= k\pi$, then $s(\tau)=\frac{2}{1-\tau}$.

\noindent
ii)If $d= k\pi$, then $s(\tau)=\frac{2}{5-\tau}$.\\

In {\bf\cite{1}}, Kerdache et al. investigated the energy decay of a coupled system of wave equations
subject to two fractional boundary damping mechanisms. More precisely, the system is given by
$$
\left\{
\matrix{
u_{tt}(x,t)-u_{xx} (x,t)+\alpha (u-v)=0 \hfill &\hbox{ in }
(0, 1)\times (0,+\infty ), \hfill & \cr
v_{tt}(x,t)-v_{xx}+\alpha (v-u)=0\hfill & \hbox{ in }%
(0, 1)\times (0,+\infty ),\hfill & \cr
u(0,t)=v(0,t)=0 \hfill & \hbox{ on } (0,+\infty ),\hfill & \cr
u_{x}(1,t)+\rho \partial _{t}^{\tau ,\omega }u(1,t)=0
\hfill & \hbox{ on } (0,+\infty ),\hfill & \cr
v_{x}(1,t)+\tilde{\rho}\partial _{t}^{\tau ,\omega
}u(1,t)=0\hfill & \hbox{ on }%
(0,1)\times (0,+\infty ),\hfill & \cr
u(x,0)=u_{0}(x),\ u_{t}(x,0)=u_{1}(x)\hfill & \hbox{ on }
(0,1),\hfill & \cr
v(x,0)=v_{0}(x),\ v_{t}(x,0)=v_{1}(x)\hfill & \hbox{ on }(0,1).\hfill & \cr}
\right.
$$
By means of semigroup theory, they established an optimal polynomial energy decay rate.

Very recently, Benzaid and Benaissa {\bf\cite{benz}} investigated the well-posedness and indirect
boundary stabilization of a one-dimensional system of weakly coupled degenerate wave equations subject to a single fractional boundary damping
acts at nondegenerate boundary endpoint, that is,
$$
\left\{
\matrix{
u_{tt}(x,t)-(a(x)u_{x})_{x}(x,t)+\alpha v=0 \hfill & \hbox{ in }(0,1)\times (0,+\infty ), \hfill & \cr
v_{tt}(x,t)-(a(x)v_{x})_{x}(x,t)+\alpha u=0\hfill & \hbox{ in }(0,1)\times (0,+\infty ),\hfill & \cr
\left\{\matrix{u(0, t)=0 \hfill & \hbox{ if } 0\leq m_a< 1 \hfill & \cr
(a u_{x})(0,t)=0 \hfill & \hbox{ if } 1\leq m_a< 2 \hfill & \cr}\right.\hfill  & \hbox{ in } (0, +\infty),\hfill&\cr
\left\{\matrix{v(0, t)=0 \hfill & \hbox{ if } 0\leq m_a< 1 \hfill & \cr
(a v_{x})(0,t)=0 \hfill & \hbox{ if } 1\leq m_a< 2 \hfill & \cr}\right.\hfill  & \hbox{ in } (0, +\infty),\hfill&\cr
v(1, t)=0 \hfill & \hbox{ for } t\in (0, +\infty),\hfill & \cr
\beta u(1,t)+ (a u_{x})(1, t)=- \rho \partial_{t}^{\tau, \omega} u(1, t) \hfill & \hbox{ in } (0, +\infty),\hfill&\cr
u(x,0)=u_0(x), u_{t}(x,0)=u_1(x)\hfill & \hbox{ for } x\in (0, 1),\hfill &\cr
v(x,0)=v_0(x), v_{t}(x,0)=v_1(x)\hfill & \hbox{ for } x\in (0, 1).\hfill &\cr}
\right.
\leqno{(ND)}
$$
The authors established the well-posedness of the system and proved its strong asymptotic stability. Furthermore,
they obtained an optimal polynomial energy decay rate of the order $t^{-\frac{2}{3-\tau}}$
by employing a frequency-domain approach combined with suitable multiplier techniques. Moreover, the degeneracy does not
affect the decay rates of the energy.

We emphasize that the multiplier method, which is commonly used in the study of boundary stabilization problems, is not applicable when the damping
mechanism acts at the degenerate endpoint. Moreover, an important distinction between the work
in {\bf\cite{zerk}} and the present study is that the former deals with a single degenerate wave equation, whereas we consider a weakly coupled system
of degenerate wave equations subject to a single fractional boundary damping mechanism. In this setting, the indirect transfer of dissipation
from the damped component to the undamped one constitutes a major analytical challenge.

Motivated by the aforementioned works {\bf\cite{lira}}, {\bf\cite{beai}},  {\bf\cite{zerk}} and {\bf\cite{benz}}, we aim to investigate the asymptotic behavior
of weakly coupled degenerate wave equations in the presence of a single fractional boundary damping mechanism acting on only one of the equations and at the degenerate endpoint.
Our main objective is to establish the stability properties of the coupled system and to characterize the long-time behavior of its energy.

The remainder of this paper is organized as follows. Section 2 introduces the appropriate functional spaces naturally associated with degenerate problems,
together with some preliminary results that will be used throughout the paper. Section 3 is devoted to establishing the well-posedness of the considered
system and investigating its strong stability properties. Section 4 is concerned with establishing the absence of exponential stability.
Finally, in Sections 5 and 6, we prove the optimal polynomial energy decay rate $t^{-\frac{2}{2+2\nu_{\gamma}-\tau}}$ for sufficiently regular initial data,
by combining the frequency-domain approach with a Bessel-function representation of the resolvent operator.

\section{Preliminary results}
We now introduce the following weighted Sobolev spaces, following the framework developed in {\bf\cite{canna.5}}:
$$
H_{0, \gamma}^{1}(0,1)=\left\{u \hbox{ is locally absolutely
continuous in } (0, 1]: x^{\gamma/2}u_x \in L^2(0,1)/\
u(1)=0\right\},
$$
$$
H_{00, \gamma}^{1}(0,1)=\left\{u\in H_{\gamma}^{1}(0, 1) /\ u(0)=u(1)=0\right\}.
$$
$$
H_{\gamma}^{1}(0,1)=\left\{u \hbox{ is locally absolutely continuous in } (0, 1]: x^{\gamma/2}u_x \in L^2(0,1)\right\}.
$$
We note that $H_\gamma^1(0,1)$ is a Hilbert space when endowed with the inner product
$$
(u, v)_{H_\gamma^1(0, 1)}=\Int_{0}^{1}(u\overline{v}+ x^{\gamma} u_x(x)\overline{v_x(x)})\, dx,\quad\forall u, v\in H_\gamma^1(0, 1).
$$
We also define
$$
|u|_{H_{0, \gamma}^1(0, 1)}=\left(\Int_{0}^{1}x^{\gamma} |u_x(x)|^2\,
dx\right)^{1/2}\quad \forall u\in H_\gamma^1(0, 1).
$$
In fact, $|\cdot|_{H_{0, \gamma}^1(0, 1)}$ defines a norm on $H_{0,\gamma}^1(0,1)$ that is equivalent to the norm inherited from $H_\gamma^1(0,1)$.
This equivalence follows directly from the following version of Poincar\'e's inequality.
\begin{proposition} There is a positive constant $C_*=C(\gamma)$ such that
\begin{equation}
\|u\|_{L^{2}(\Omega)}^2\leq C_* |u|_{H_{0, \gamma}^1(0, 1)}^2\quad
\forall u\in H_{0, \gamma}^1(0, 1). \label{eh10}
\end{equation}
\end{proposition}
{\bf Proof.} Let $u\in H_{0, \gamma}^1(0, 1)$. For any $x \in ]0, 1]$ we
have that
$$
|u(x)|=\left|\Int_{x}^{1}u_x(s)\, ds\right|\leq |u|_{H_{0, \gamma}^1(0, 1)}\left\{\Int_{0}^{1}\Frac{1}{x^\gamma}\,
ds\right\}^{1/2}.
$$
Therefore
$$
\Int_{0}^{1}|u(x)|^2\, dx\leq \Frac{1}{1-\gamma} |u|_{H_{0, \gamma}^1(0, 1)}^2.
$$
Next, we define
$$
H_{\gamma}^2(0, 1)=\{u\in H_{\gamma}^1(0, 1):\ x^\gamma u_x\in H^1 (0, 1)\},
$$
where $H^1 (0, 1)$ denotes the classical Sobolev space.
\begin{remarque}
We observe that, for $u\in H_\gamma^2(0,1)$ with $\gamma\in[1,2)$, the weighted flux necessarily satisfies
$(x^\gamma u_x)(0)=0$, since $x^{-\gamma}$ is not integrable in any neighbourhood of the origin. Consequently, the boundary condition at $x=0$ cannot be
prescribed independently, and the corresponding problem cannot be formulated as a well-posed semigroup evolution problem on the Hilbert space under consideration.
\label{rem1}
\end{remarque}
The elements of $H_{0, \gamma}^1(0, 1)$ satisfy the following property.
\begin{proposition}
For every $u\in H_{0, \gamma}^{1}(0, 1)$, $u$ is absolutely continuous in $[0, 1]$ and we have
$$
\|u\|_{L^{\infty}(0, 1)}\leq \frac{1}{\sqrt{1-\gamma}}\|u\|_{H_{0, \gamma}^{1}(0, 1)}.
\leqno{(*)}
$$
\label{ppa1}
\end{proposition}
{\bf Proof.} As
$$
u'(x)=\Frac{1}{x^{\gamma/2}} x^{\gamma/2}u'(x)\qquad \forall x\in ]0, 1].
$$
then
$$
\matrix{\Int_{0}^{1}|u'(x)|\, dx
&\leq &\left(\Int_{0}^{1}\Frac{1}{x^{\gamma}}\, dx\right)^{1/2} \|u\|_{H_{0, \gamma}^{1}(0, 1)}  \hfill \cr
&=& \Frac{1}{\sqrt{1-\gamma}}\|u(x)\|_{H_{0, \gamma}^{1}(0, 1)}\hfill \cr}
$$
$u'$ is summable over $(0, 1)$. So $u$ is absolutely continuous in $[0, 1]$. Moreover
$$
|u(x)|=\left|\Int_{x}^{1}u'(s)\, ds\right|\leq \Int_{0}^{1}|u'(x)|\, dx\leq \Frac{1}{\sqrt{1-\gamma}}\|u\|_{H_{0, \gamma}^{1}(0, 1)}.
$$
$(*)$ is thus proved.

\section{Well-posedness and strong stability}
\subsection{Augmented model}
In this section we reformulate (\ref{eqq100}) into an augmented system. For that, we need the following proposition.
\begin{proposition}[see {\bf\cite{1}}]
Let $\vartheta$ be the function:
\begin{equation}
\vartheta(\varsigma)=|\varsigma|^{(2\tau-1)/2},\quad -\infty<\varsigma<+\infty,\
0<\tau<1. \label{e1}
\end{equation}
Then the relationship between the `input' U and the `output' O of
the system
\begin{equation}
\partial_t\varphi(\varsigma, t)+(\varsigma^{2}+\omega)\varphi(\varsigma, t) -U(t)\vartheta(\varsigma)=0,\quad -\infty<\varsigma<+\infty,\omega\geq 0, t> 0,
\label{e6}
\end{equation}
\begin{equation}
\varphi(\varsigma, 0)=0, \label{e7}
\end{equation}
\begin{equation}
O(t)=(\pi)^{-1}\sin(\tau\pi)\Int_{-\infty}^{+\infty}\vartheta(\varsigma)\varphi(\varsigma, t)\, d\varsigma,
\label{e8}
\end{equation}
where $U\in C^0([0,+\infty))$, is given by
\begin{equation}
O=I^{1-\tau, \omega}U, \label{e9}
\end{equation}
where
$$
[I^{\tau,
\omega}f](t)=\Frac{1}{\Gamma(\tau)}\Int_{0}^{t}(t-s)^{\tau-1}e^{-\omega(t-s)}f(s)\, ds.
$$
\label{th2}
\end{proposition}
\begin{lemma}[see {\bf\cite{1}}]
If $\lambda\in D_{\omega}=\C\backslash]-\infty, -\omega]$ then
$$
\Int_{-\infty}^{+\infty}\Frac{\vartheta^2(\varsigma)}{\lambda+\omega+\varsigma^2}\, d\varsigma
=\Frac{\pi}{\sin\tau\pi}(\lambda+\omega)^{\tau-1}.
$$
\label{achour}
\end{lemma}
We are now ready to reformulate system (\ref{eqq100}). Indeed, in view of Proposition \ref{th2}, system (\ref{eqq100})
can be equivalently rewritten as the following augmented system:
$$
\left\{
\matrix{
u_{tt}(x,t)-(x^{\gamma}u_{x})_{x}(x,t)+\alpha v=0 \hfill & \hbox{ in }(0,1)\times
(0,+\infty ), \hfill & \cr
v_{tt}(x,t)-(x^{\gamma}v_{x})_{x}(x,t)+\alpha u=0\hfill & \hbox{ in }(0,1)\times
(0,+\infty ),\hfill & \cr
\varphi_t(\varsigma, t)+(\varsigma^{2}+\omega)\varphi(\varsigma, t) -u_{t}(0, t)\vartheta(\varsigma)=0, \hfill & \varsigma\in \R,\hfill \omega\geq 0,\hfill t>0,\hfill&\cr
(x^{\gamma}u_{x})(0,t)=\zeta \Int_{-\infty}^{+\infty}\vartheta(\varsigma)\varphi(\varsigma, t)\, d\varsigma \hfill  & \zeta=\rho(\pi)^{-1}\sin(\tau\pi) ,\hfill&\cr
u(1, t)=v(0, t)=v(1, t)=0 \hfill & \hbox{ for } t\in (0, +\infty),\hfill & \cr
u(x,0)=u_0(x), u_{t}(x,0)=u_1(x)\hfill & \hbox{ for } x\in (0, 1),\hfill &\cr
v(x,0)=v_0(x), v_{t}(x,0)=v_1(x)\hfill & \hbox{ for } x\in (0, 1).\hfill &\cr}
\right.
\leqno{(P')}
$$
The energy functional associated with solutions of problem $(P')$ is defined by
\begin{equation}
\matrix{{\cal E}(t)=\Frac{1}{2}\int_{0}^{1}(|u_t|^2+x^{\gamma} |u_x|^2)dx+\Frac{1}{2}\int_{0}^{1}(|v_t|^2+x^{\gamma} |v_x|^2)dx\hfill &\cr
\Frac{1}{2}\alpha\int_{0}^{1}(u {\overline v}+v{\overline u})dx
+\Frac{\zeta}{2}\int_{-\infty}^{+\infty}|\varphi(\varsigma, t)|^{2}\, d\varsigma. & \cr}
\label{e10}
\end{equation}
\begin{lemma}
Let $(u, v, \varphi)$ be a regular solution of the problem $(P')$.
Then, the energy functional defined by (\ref{e10}) satisfies
\begin{equation}
{\cal E}'(t)=-\zeta \int_{-\infty}^{+\infty}(\varsigma^{2}+\omega)|\varphi(\varsigma, t)|^{2}\, d\varsigma\leq 0.
\label{e11}
\end{equation}
\label{lem1}
\end{lemma}
In this section, we establish the existence and uniqueness of solutions to problem $(P')$ within the
framework of semigroup theory. To this end, we introduce the state vector
$U=(u, {\tilde u}, v, {\tilde v}, \varphi)^{T}$, where ${\tilde u}=u_t$, ${\tilde v}=v_t$,
Accordingly, system $(P')$ can be reformulated as an abstract Cauchy evolution problem.
\begin{equation}
\left\{\matrix{U'={\cal P}U,\qquad\hbox{ for all } t> 0, \hfill & \cr
U(0)=U_0,\hfill & \cr}
\right.
\label{e14}
\end{equation}
where $U_0=(u_0, u_1, v_0, v_1, \varphi_{0})^{T}$ and
$$
{\cal P}: D({\cal P})\subset {\cal H}\rightarrow {\cal H}
$$
is the operator given by
\begin{equation}
{\cal P}\pmatrix{u\cr {\tilde u}\cr v\cr {\tilde v}\cr \varphi\cr}=
\pmatrix{{\tilde u}\cr (x^{\gamma} u_x)_x-\alpha v\cr {\tilde v}\cr (x^{\gamma}v_x)_x-\alpha u\cr
-(\varsigma^{2}+\omega)\varphi+{\tilde u}(0)\vartheta(\varsigma)\cr}.
\label{e1972}
\end{equation}
We introduce the following phase space, which will serve as the natural energy space for the system:
$$
{\cal H}=H_{0, \gamma}^{1}(0,1)\times L^{2}(0,1)\times H_{00, \gamma}^{1}(0,1)\times L^{2}(0,1)\times L^{2}(-\infty, +\infty).
$$
Endowed with the following inner product, ${\cal H}$ is a Hilbert space.
$$
\matrix{\langle U, {\tilde U} \rangle_{\cal H}=\Int_{0}^{1}x^{\gamma} u_{1x} {\overline u}_{2x} dx+\Int_{0}^{1}x^{\gamma}v_{1x} {\overline v}_{2x} dx
+\alpha\Int_{0}^{1} (u_1 {\overline v}_2+v_2 {\overline u}_1)\, dx\hfill & \cr
+\Int_{0}^{1}{\tilde u}_1 {\overline {\tilde u}}_2 dx +\Int_{0}^{1}{\tilde v}_1 {\overline {\tilde v}}_2 dx
+\zeta\Int_{-\infty}^{+\infty} \varphi_1  {\overline \varphi_2}  \, d\varsigma,\hfill &\cr}
$$
for all $U= (u_1, {\tilde u}_1, v_1, {\tilde v}_1, \varphi_1)^{T}$ and $\tilde{U}=(u_2, {\tilde u}_2, v_2, {\tilde v}_2, \varphi_2)^{T}$.

The domain of ${\cal P}$ is
\begin{equation}
D({\cal P})=\left\{\matrix{(u, {\tilde u}, v,{\tilde v}, \varphi)^{T}\ \hbox{in } {\cal H}:
u\in H_{\gamma}^{2}(0,1)\cap H_{0, \gamma}^{1}(0,1), v\in H_{\gamma}^{2}(0,1)\cap H_{00, \gamma}^{1}(0,1),  \hfill \cr
{\tilde u}\in H_{0, \gamma}^{1}(0,1), {\tilde v}\in H_{00, \gamma}^{1}(0,1),
-(\varsigma^{2}+\omega)\varphi+{\tilde u}(0)\vartheta(\varsigma)\in L^{2}(-\infty, +\infty)  ,\hfill \cr
(x^{\gamma} u_x)(0)
=\zeta\Int_{-\infty}^{+\infty}\vartheta(\varsigma)\varphi(\varsigma)\, d\varsigma, \hfill \cr
|\varsigma|\varphi\in L^{2}(-\infty, +\infty)\hfill \cr
}\right\}.
\label{e15}
\end{equation}
We are now in a position to state the following existence and uniqueness result.
\begin{theorem}[Existence and uniqueness]\
\begin{itemize}
\item[(1)] If $U_0\in D({\cal P})$, then system (\ref{e14}) has a unique strong solution with the following regularity,
$$
U\in C^0(\R_{+}, D({\cal P}))\cap C^1(\R_{+}, {\cal H}).
$$
\item[(2)] If $U_0\in {\cal H}$, then system (\ref{e14}) has a unique weak solution such that
$$
U\in C^0(\R_{+}, {\cal H}).
$$
\end{itemize}
\label{1th}
\end{theorem}
{\bf Proof.}\\
We adopt the semigroup approach. To this end, we first show that the operator ${\cal P}$ is monotone.
Let $U\in D({\cal P})$. By virtue of (\ref{e14}) and (\ref{e11}), together with the fact that
\begin{equation}
{\cal E}(t)=\Frac{1}{2}\|U\|_{\cal H}^{2},
\label{e16}
\end{equation}
we have
\begin{equation}
\Re\langle {\cal P} U,  U \rangle_{\cal H}=-\zeta \Int_{-\infty}^{+\infty}(\varsigma^{2}+\omega)|\varphi(\varsigma)|^{2}\, d\varsigma.
\label{e17}
\end{equation}
Consequently, the operator ${\cal P}$ is dissipative.
We next establish the surjectivity of $(\lambda I-{\cal P})$ for every $\lambda>0$. More precisely, given
\begin{equation}
(\lambda I-{\cal P})U={\cal G}.
\label{e18z}
\end{equation}
From Equation (\ref{e18z}), we get the following system of equations
\begin{equation}
\left\{\matrix{\lambda u-{\tilde u}=g_1,  \hfill & \cr
\lambda {\tilde u}-(x^{\gamma}u_x)_x+\alpha v=g_2, \hfill & \cr
\lambda v-{\tilde v}=g_3\hfill & \cr
\lambda {\tilde v}-(x^{\gamma}v_x)_x+\alpha u=g_4, \hfill & \cr
\lambda \varphi+(\varsigma^{2}+\omega)\varphi-{\tilde u}(0)\vartheta(\varsigma)=g_5.\hfill & \cr}\right.
\label{e18}
\end{equation}
Assume that $u$ and $v$ have been constructed with the required regularity. Then, $(\ref{e18})_1$ and $(\ref{e18})_3$ yield
\begin{equation}
\left\{\matrix{{\tilde u}=\lambda u-g_1 \in H_{0, \gamma}^{1}(0,1),  \hfill & \cr
{\tilde v}=\lambda v-g_3 \in H_{00, \gamma}^{1}(0,1),  \hfill & \cr}\right.
\label{e19}
\end{equation}
By using $(\ref{e18})_2, (\ref{e18})_4$ and (\ref{e19}), we readily obtain that $u$ and $v$ satisfy
\begin{equation}
\left\{\matrix{\lambda^{2}u -(x^{\gamma}u_x)_x+\alpha v=g_2+\lambda g_1,\hfill & \cr
\lambda^{2}u -(x^{\gamma}u_x)_x+\alpha u=g_4+\lambda g_3.\hfill & \cr}\right.
\label{e21}
\end{equation}
Solving system (\ref{e21}) is equivalent to finding a pair $(u,v)$ in $H_{\gamma}^{2}(0,1)\cap H_{0, \gamma}^{1}(0,1)\times
H_{\gamma}^{2}(0,1)\cap H_{00, \gamma}^{1}(0,1)$ such that
\begin{equation}
\left\{\matrix{\Int_{0}^{1}(\lambda^{2}u {\overline w}- (x^{\gamma}u_x)_x {\overline w})\, dx+\alpha \Int_{0}^{1}v{\overline w}\, dx
=\Int_{0}^{1}(g_2+\lambda g_1){\overline w}\, dx,\hfill & \cr
\Int_{0}^{1}(\lambda^{2}v {\overline y}- (x^{\gamma} v_x)_x {\overline y})\, dx+\alpha \Int_{0}^{1}u{\overline y}\, dx
=\Int_{0}^{1}(g_4+\lambda g_3){\overline y}\, dx,\hfill & \cr}\right.
\label{e22}
\end{equation}
for all $(w, y)\in H_{0, \gamma}^{1}(0,1)\times H_{00, \gamma}^{1}(0,1)$.
By combining (\ref{e22}) with the boundary condition $(\ref{e15})_3$ and
relation $(\ref{e18})_5$, we obtain the following system satisfied by the functions $u$ and $v$.
\begin{equation}
\left\{\matrix{
\matrix{ \Int_{0}^{1}(\lambda^{2}u {\overline w}+ x^{\gamma}u_x
{\overline w_x})\, dx +\alpha \Int_{0}^{1}v{\overline w}\, dx+{\tilde\zeta} {\tilde u}(0) {\overline w(0)}\, \hfill &\cr
\quad\ \qquad=\Int_{0}^{1}(g_2+\lambda g_1){\overline w}\, dx
-\zeta\Int_{-\infty}^{+\infty} \Frac{\vartheta(\varsigma)}{\varsigma^{2}+\omega+\lambda}g_5(\varsigma)\,  d\varsigma{\overline w(0)}, &\cr}\hfill & \cr
\matrix{ \Int_{0}^{1}(\lambda^{2}v {\overline y}+ x^{\gamma}v_x
{\overline y_x})\, dx+\alpha \Int_{0}^{1}u{\overline y}\, dx=\Int_{0}^{1}(g_4+\lambda g_3){\overline y}\, dx, &\cr}\hfill & \cr}\right.
\label{e23}
\end{equation}
where ${\tilde\zeta}=\zeta\Int_{-\infty}^{+\infty}
\Frac{\vartheta^{2}(\varsigma)}{\varsigma^{2}+\omega+\lambda}  \, d\varsigma$. Using again
$(\ref{e19})_1$, we deduce that
\begin{equation}
{\tilde u}(0)=\lambda u(0)-g_1(0). \label{e24}
\end{equation}
Inserting (\ref{e24}) into (\ref{e23}), we get
\begin{equation}
\left\{\matrix{
\matrix{ \Int_{0}^{1}(\lambda^{2}u {\overline w}+  x^{\gamma} u_x
{\overline w_x})\, dx+\alpha \Int_{0}^{1}v{\overline w}\, dx +\lambda {\tilde\zeta} u(0){\overline w(0)}\, \hfill &\cr
=\Int_{0}^{1}(g_2+\lambda g_1){\overline w}\, dx-\zeta \Int_{-\infty}^{+\infty}
\Frac{\vartheta(\varsigma)}{\varsigma^{2}+\omega+\lambda}g_5( \varsigma)\,  d\varsigma{\overline w(0)} +{\tilde\zeta} g_1(0) {\overline w(0)}, &\cr}\hfill &\cr
\Int_{0}^{1}(\lambda^{2}v {\overline y}+ x^{\gamma}v_x
{\overline y_x})\, dx+\alpha \Int_{0}^{1}u{\overline y}\, dx=\Int_{0}^{1}(g_4+\lambda g_3){\overline y}\, dx.
\hfill & \cr}\right.
\label{e25}
\end{equation}
Adding $(\ref{e25})_1$ and $(\ref{e25})_2$, we introduce a sesquilinear form ${\cal B}: [H_{0, \gamma}^{1}(0,1)\times H_{00, \gamma}^{1}(0,1)]^2\rightarrow
\C$ given by
$$
{\cal B}((u, v), (w, y))=\lambda^{2}\Int_{0}^{1}(u {\overline w}+v {\overline y})\, dx
+ \Int_{0}^{1}x^{\gamma}(u_x{\overline w_x}+v_x{\overline y_x})\, dx
+\alpha \Int_{0}^{1}(v{\overline w}+u{\overline y})\, dx +\lambda {\tilde\zeta} u(0){\overline w(0)},
$$
and a continuous antilinear functional ${\cal L}: (0,1) H_{0, \gamma}^{1}(0,1)\times H_{00, \gamma}^{1}(0,1)\rightarrow \C$ where
$$
{\cal L}(w, y)= \Int_{0}^{1}(g_2+\lambda g_1){\overline w}\, dx-\zeta
\Int_{-\infty}^{+\infty} \Frac{\vartheta(\varsigma)}{\varsigma^{2}+\omega+\lambda}g_5(\varsigma)\,  d\varsigma{\overline w(0)}
+{\tilde\zeta} g_1(0) {\overline w(1)}+\Int_{0}^{1}(g_4+\lambda g_3){\overline y}\, dx,
$$
satisfying
\begin{equation}
{\cal B}((u, v), (w, y))={\cal L}(w, y).
\label{e26}
\end{equation}
The sesquilinear form ${\cal B}(., .)$ is a bounded since for all
$(u, v), (w, y)\in H_{0, \gamma}^{1}(0,1)\times H_{00, \gamma}^{1}(0,1)$, we have
$$
\matrix{{\cal B}((u, v), (w, y))\leq \lambda^2 \|u\|_{L^2(0, 1)}\|w\|_{L^2(0, 1)}+\lambda^2 \|v\|_{L^2(0, 1)}\|y\|_{L^2(0, 1)}\hfill & \cr
+\|x^{\gamma/2}u_x\|_{L^2(0, 1)}\|x^{\gamma/2}w_x\|_{L^2(0, 1)}
+\|x^{\gamma/2}v_x\|_{L^2(0, 1)}\|x^{\gamma/2}y_x\|_{L^2(0, 1)}\hfill & \cr
+|\alpha|\|v\|_{L^2(0, 1)}\|w\|_{L^2(0, 1)}+|\alpha|\|u\|_{L^2(0, 1)}\|y\|_{L^2(0, 1)}
+\lambda {\tilde\zeta}|u(0)|^2\hfill & \cr
\leq M \|(u, v)\|_{H_{0, \gamma}^{1}\times H_{00, \gamma}^{1}(0,1)}\|(w, y)\|_{H_{0, \gamma}^{1}\times H_{00, \gamma}^{1}(0,1)}\hfill & \cr}
$$
and is coercive because $\forall (u, v)\in H_{0, \gamma}^{1}(0,1)\times H_{00, \gamma}^{1}(0,1)$, we have
$$
\matrix{{\cal B}((u, v), (u, v))&\geq &\lambda^2 (\|u\|_{L^2(0, 1)}^2+\|v\|_{L^2(0, 1)}^2)+
\|x^{\gamma/2}u_x\|_{L^2(0, 1)}^2+\|x^{\gamma/2}v_x\|_{L^2(0, 1)}^2\hfill  \cr
&& +2\alpha\Re\Int_{0}^{1}u{\overline v}\, dx+\lambda {\tilde\zeta} |u(0)|^2\hfill  \cr
&\geq & c (\|u\|_{H_{00, \gamma}^{1}(0,1)}+\|v\|_{H_{00, a}^{1}(0,1)}),\hfill  \cr}
$$
for $\alpha$ small enough.
Therefore, by the Lax-Milgram Theorem, system (\ref{e26}) admits a unique solution $(u, v)\in H_{0, \gamma}^{1}(0,1)\times H_{00, \gamma}^{1}(0,1)$.

Choosing $(w, y)=(w, 0)$ with $w\in {\cal D}(0, 1)$ in (\ref{e26}), we obtain
\begin{equation}
\lambda^{2}u -(x^{\gamma}u_x)_x+\alpha v=g_2+\lambda g_1.
\label{mi1}
\end{equation}
Since $u\in H_{0, \gamma}^{1}(0,1)$ and $v\in H_{00, \gamma}^{1}(0,1)$, it follows that $(x^{\gamma}u_x)_x\in L^{2}(0, 1)$. Consequently,
$u\in H_{\gamma}^{2}(0,1)\cap H_{0, \gamma}^{1}(0,1)$.

Similarly, Choosing $(w,y)=(0,y)$, with $y\in {\cal D}(0,1)$, in (\ref{e26}), yields
\begin{equation}
\lambda^{2}v -(x^{\gamma}v_x)_x+\alpha u=g_4+\lambda g_3
\label{mi2}
\end{equation}
Hence $v\in H_{\gamma}^{2}(0,1)\cap H_{00, \gamma}^{1}(0,1)$.

Next, multiplying the complex conjugates of (\ref{mi1}) and (\ref{mi2}) by $w\in H_{0,\gamma}^{1}(0,1)$
and $y\in H_{00,\gamma}^{1}(0,1)$, respectively, and integrating by parts over $(0,1)$, we obtain, upon comparison with (\ref{e26}),
$$
\matrix{ (x^{\gamma}u_x)(0){\overline w}(0)
+\rho\lambda(\lambda+\omega)^{\tau-1} u(0){\overline w}(0)\hfill & \cr
+\zeta\Int_{-\infty}^{+\infty} \Frac{\vartheta(\varsigma)}{\varsigma^{2}+\omega+\lambda}g_5(\varsigma)  \, d\varsigma{\overline w}(0)\
-\rho(\lambda+\omega)^{\tau-1}  g_1(0){\overline w}(0)=0.\hfill & \cr}
$$
Therefore, setting ${\tilde u}= \lambda u- g_1$ and recalling the definition of $\varphi$ in $(\ref{e18})_5$, we deduce that
$$
(x^{\gamma}u_x)(0)=\zeta\Int_{-\infty}^{+\infty}\vartheta(\varsigma)\varphi(\varsigma)\, d\varsigma.
$$
To complete the proof that $U\in D({\cal P})$, it remains to establish that $\varphi$ and $|\varsigma|\varphi\in L^{2}(-\infty, \infty)$.
By virtue of $(\ref{e18})_5$, we have
$$
\Int_{\R}|\varphi(\varsigma)|^2\, d\varsigma\leq 3\Int_{\R}\Frac{|g_5(\varsigma)|^2}{(\varsigma^{2}+\omega+\lambda)^2}  \, d\varsigma
+3 (\lambda^2|u(0)|^2+|g_1(0)|^2)\Int_{\R}\Frac{|\varsigma|^{2\tau-1}}{(\varsigma^{2}+\omega+\lambda)^2}  \, d\varsigma.
$$
By Proposition \ref{th2} the latter integral can be explicitly evaluated as
$$
\Int_{\R}\Frac{|\varsigma|^{2\tau-1}}{(\varsigma^{2}+\omega+\lambda)^2}  \, d\varsigma=(1-\tau)\Frac{\pi}{\sin \tau\pi}(\lambda+\omega)^{\tau-2}.
$$
On the other hand, since $g_5\in L^2(\R)$, we have
$$
\Int_{\R}\Frac{|g_5(\varsigma)|^2}{(\varsigma^{2}+\omega+\lambda)^2}  \, d\varsigma\leq \Frac{1}{(\omega+\lambda)^2}
\Int_{\R}|g_5(\varsigma)|^2  \, d\varsigma< +\infty.
$$
It follows that $\varphi\in L^{2}(\R)$. Again, using $(\ref{e18})_5$, we obtain
$$
\Int_{\R}|\varsigma\varphi(\varsigma)|^2\, d\varsigma\leq 3\Int_{\R}\Frac{|\varsigma|^{2}|g_5(\varsigma)|^2}{(\varsigma^{2}+\omega+\lambda)^2}  \, d\varsigma
+3 (\lambda^2|u(0)|^2+|g_1(0)|^2)\Int_{\R}\Frac{|\varsigma|^{2\tau+1}}{(\varsigma^{2}+\omega+\lambda)^2}  \, d\varsigma.
$$
Applying Proposition \ref{th2}once more, we find that
$$
\Int_{\R}\Frac{|\varsigma|^{2\tau+1}}{(\varsigma^{2}+\omega+\lambda)^2}  \, d\varsigma=\tau\Frac{\pi}{\sin \tau\pi}(\lambda+\omega)^{\tau-1}.
$$
Furthermore, since $g_5\in L^2(\R)$, it follows that
$$
\Int_{\R}\Frac{|\varsigma|^{2}|g_5(\varsigma)|^2}{(\varsigma^{2}+\omega+\lambda)^2}  \, d\varsigma\leq \Frac{1}{(\omega+\lambda)}
\Int_{\R}|g_5(\varsigma)|^2  \, d\varsigma< +\infty.
$$
Consequently, $|\varsigma|\varphi\in L^{2}(\R)$. Finally, recalling that $\varphi\in L^{2}(\R)$, we obtain
$$
-(\varsigma^{2}+\omega)\varphi+{\tilde u}(0)\vartheta(\varsigma)=\lambda\varphi(\varsigma)-g_5(\varsigma)\in L^{2}(\R).
$$
Thus, all the required regularity and boundary conditions are satisfied, and hence $U\in D({\cal P})$.
We conclude that the operator $\lambda I-{\cal P}$ is surjective for every $\lambda>0$.

\hfill$\diamondsuit$\\

\subsection{Strong stability of the system}
n this section, we apply the general criterion established in Theorem \ref{thm14} to prove the strong stability
of the $C_0$-semigroup $e^{t{\cal P}}$ generated by the operator ${\cal P}$ associated with the wave system (\ref{eqq100})
in the absence of the compactness of the resolvent of ${\cal P}$.

To state and establish our stability results, we first recall some preliminary results from semigroup theory.
\begin{theorem}[\cite{arba}]
Let ${\cal P}$ be the generator of a uniformly bounded
$C_0$-semigroup $\{S(t)\}_{t\geq 0}$ on a Hilbert space ${\cal X}$.
If:
\begin{itemize}
\item[(i)] ${\cal P}$ does not have eigenvalues on $i\R$.
\item[(ii)] The intersection of the spectrum $\sigma({\cal P})$ with $i\R$ is at most a countable set,
\end{itemize}
then the semigroup $\{S(t)\}_{t\geq 0}$ is asymptotically stable,
i.e, $\|S(t)z\|_{{\cal X}}\rightarrow 0$ as $t\rightarrow \infty$
for any $z\in {\cal X}$.
\label{thm14}
\end{theorem}
We are now in a position to state our main result.
\begin{theorem}
The $C_0$-semigroup $e^{t{\cal P}}$ is strongly stable in $\mathcal{H}$ if and only if the coefficient $\alpha$
satisfies
$$
\alpha\not=\Frac{1}{2}\left(\Frac{2-\gamma}{2}\right)^2 (j_{\nu_{\gamma}, k}^2-j_{\nu_{\gamma}, m}^2),\quad
k, m\in \N,
\leqno{(C)}
$$
where $\nu_{\gamma}=(1-\gamma)/(2-\gamma)$ and $j_{\nu, 1}< j_{\nu, 2}< \ldots <j_{\nu, k}<\ldots$ denote the sequence of positive zeros
of the Bessel function of first kind and of order $\nu$.
\label{thmba1}
\end{theorem}
For the proof of Theorem \ref{thmba1}, we need the following two lemmas.
\begin{lemma}
${\cal P}$ does not have eigenvalues on $i\R$.
\label{lem10}
\end{lemma}
{\bf Proof.}
Let $\lambda\in \R-\{0\}$. We prove that $i\lambda$ is not an eigenvalue of ${\cal P}$.
Let $U=(u, {\tilde u}, v, {\tilde v}, \varphi)^{T}$ with $\|U\|_{{\cal H}}=1$, be such that
\begin{equation}
{\cal P}U=i\lambda U.
\label{wwzz}
\end{equation}
Using the definition of ${\cal P}$ it follows that ${\cal P}U=i\lambda U$ if and only if
\begin{equation}
\left\{\matrix{i\lambda u-{\tilde u}=0,  \hfill & \cr
i\lambda {\tilde u}-(x^{\gamma}u_x)_x+\alpha v=0, \hfill & \cr
i\lambda v-{\tilde v}=0\hfill & \cr
i\lambda {\tilde v}-(x^{\gamma}v_x)_x+\alpha u=0, \hfill & \cr
i\lambda \varphi+(\varsigma^{2}+\omega)\varphi-{\tilde u}(0)\vartheta(\varsigma)=0.\hfill & \cr}\right.
\label{eqq96}
\end{equation}
Thus, by eliminating ${\tilde u}, {\tilde v}$, the system (\ref{eqq96}) implies that
\begin{equation}
\left\{\matrix{
\lambda^2 u+(x^{\gamma}u_x)_x-\alpha v=0, \hfill & \cr
\lambda^2 v+(x^{\gamma}v_x)_x-\alpha u=0, \hfill & \cr
(i\lambda+\varsigma^2+\omega)\varphi-{\tilde u}(0)\vartheta(\varsigma)=0,\hfill &\cr
v(0)=v(1)=u(1)=0, \hfill &\cr
(x^{\gamma} u_x)(0)=\zeta\Int_{-\infty}^{+\infty}\vartheta(\varsigma)\varphi(\varsigma)\, d\varsigma.\hfill & \cr
}\right.
\label{v30}
\end{equation}
Using (\ref{e17}) and (\ref{wwzz}), we find
\begin{equation}
\varphi\equiv 0,
\label{em28}
\end{equation}
then, using the third equation in $(\ref{v30})$, we deduce that
\begin{equation}
{\tilde u}(0)=0.
\label{em29}
\end{equation}
Therefore, from $(\ref{eqq96})_1$  and $(\ref{v30})_2$, we get
\begin{equation}
u(0)=0 \quad \hbox{and} \quad (x^{\gamma}u_x)(0)=0.
\label{em30}
\end{equation}
Then $\phi=u+v$ and $\psi=u-v$ satisfy
\begin{equation}
\left\{
\matrix{
(\lambda^2-\alpha)\phi+(x^{\gamma}\phi_{x})_x=0, \hfill & \cr
(\lambda^2+\alpha)\psi+(x^{\gamma}\psi_{x})_x=0. \hfill & \cr}
\right.
\label{e31qqq}
\end{equation}
The solution of the equation (\ref{e31qqq}) is given by
$$
\left\{
\matrix{
\phi(x)=c_1 \Phi_{+}(x) + c_2\Phi_{-}(x),\hfill & \cr
\psi(x)={\tilde c}_1 \Phi_{++}(x) + {\tilde c}_2\Phi_{--}(x),\hfill & \cr}
\right.
$$
togheter with the boundary conditions
$$
\phi(0)=\phi(1)=\psi(0)=\psi(1)=0, (x^{\gamma} u_x)(0)==0,
$$
where
\begin{equation}
\left\{
\matrix{
\Phi_{+}(x)=x^{\frac{1-\gamma}{2}}J_{\nu_{\gamma}}\left(\frac{2}{2-\gamma}\sqrt{\lambda^2-\alpha} x^{\frac{2-\gamma}{2}}\right),
\quad  \Phi_{-}(x)=x^{\frac{1-\gamma}{2}}J_{-\nu_{\gamma}}\left(\frac{2}{2-\gamma}\sqrt{\lambda^2-\alpha} x^{\frac{2-\gamma}{2}}\right),\hfill &\cr
\Phi_{++}(x)=x^{\frac{1-\gamma}{2}}J_{\nu_{\gamma}}\left(\frac{2}{2-\gamma}\sqrt{\lambda^2+\alpha} x^{\frac{2-\gamma}{2}}\right),\hfill &\cr
\Phi_{--}(x)=x^{\frac{1-\gamma}{2}}J_{-\nu_{\gamma}}\left(\frac{2}{2-\gamma}\sqrt{\lambda^2+\alpha} x^{\frac{2-\gamma}{2}}\right),
\hfill &\cr}
\right.
\label{mmrr}
\end{equation}
where
\begin{equation}
J_\nu(y)=\Sum_{m=0}^{\infty}\Frac{(-1)^m}{m!\Gamma(m+\nu+1)}\left(\Frac{y}{2}\right)^{2m+\nu}
=\Sum_{m=0}^{\infty}c_{\nu, m}^{+}y^{2m+\nu},
\label{bess1}
\end{equation}
\begin{equation}
J_{-\nu}(y)=\Sum_{m=0}^{\infty}\Frac{(-1)^m}{m!\Gamma(m-\nu+1)}\left(\Frac{y}{2}\right)^{2m-\nu}
=\Sum_{m=0}^{\infty}c_{\nu, m}^{-}y^{2m-\nu},
\label{bess2}
\end{equation}
where $J_{\nu}$ and $J_{-\nu}$ are Bessel functions of the first kind of order $\nu$ and $-\nu$.

As $\phi(0)=\psi(0)=0$, then $c_2={\tilde c}_2=0$. As $\phi(1)=\psi(1)=0$, then
$$
c_1 \Phi_{+}(1)=0, \ {\tilde c}_1 \Phi_{++}(1)=0.
$$
As $u(x)=\Frac{1}{2}(\phi(x)+\psi(x))$, we deduce that
$$
(1-\gamma){\tilde c}_{\nu_\gamma, 0}^{+} c_1+(1-\gamma){\widetilde{\tilde c}}_{\nu_\gamma, 0}^{+} {\tilde c}_1=0
$$
where
$$
{\tilde{\tilde c}}_{\nu_{\gamma, 0}}^{+}= c_{\nu_{\gamma, 0}}^{-}\left(\Frac{2}{2-\gamma}{\tilde \lambda}\right)^{\nu_{\gamma}},
{\tilde{\tilde c}}_{\nu_{\gamma, 0}}^{+}=c_{\nu_{\gamma, 0}}^{+}\left(\Frac{2}{2-\gamma}{\widetilde{\tilde \lambda}}\right)^{\nu_{\gamma}},
$$
where ${\tilde \lambda}=\sqrt{\lambda^2-\alpha},\  {\widetilde{\tilde \lambda}}=\sqrt{\lambda^2+\alpha}$.
If Bessel are zero then
$$
\Frac{2}{2-\gamma}\sqrt{\lambda^2-\alpha}=j_{\nu_{\gamma}, k} \hbox{ and }\Frac{2}{2-\gamma}\sqrt{\lambda^2+\alpha}=j_{\nu_{\gamma}, m}
$$
for some integers $k$ and $m$. Hence, eigenvalues on $i\R$ exist iff
$$
\alpha=\Frac{1}{2}\left(\Frac{2-\gamma}{2}\right)^2 (j_{\nu_{\gamma}, k}^2-j_{\nu_{\gamma}, m}^2).
$$
Hence, if condition $(C)$ is satisfied we deduce that $c_1 = 0$ or ${\tilde c}_1 = 0$ and consequently $u = v = 0$.

Therefore $U=0$. Consequently, ${\cal P}$ does not have purely imaginary eigenvalues.

\begin{lemma}
We have
$$
\matrix{i\R\subset \rho({\cal P}) \hbox{ if }\omega\not=0, \hfill &\cr
i\R^*\subset \rho({\cal P}) \hbox{ if }\omega=0,\hfill &\cr}
$$
\label{lem11}
\end{lemma}
where $\R^*=\R-\{0\}$.\\
\noindent
{\bf Proof.}\\
\noindent
$\bullet${\bf Case 1}: $\lambda\not=0$.\\
We will prove that the operator $i\lambda I-{\cal P}$ is surjective for $\lambda\not=0$. For this purpose, let
${\cal G}=(g_1, g_2, g_3, g_4, g_5)^{T}\in {\cal H}$, we seek $U=(u, {\tilde u}, v, {\tilde v}, \varphi)^{T}\in D({\cal P})$ solution of the following equation
\begin{equation}
(i\lambda I-{\cal P})U={\cal G}.
\label{e35eee}
\end{equation}
Equivalently, we have
\begin{equation}
\left\{\matrix{i\lambda u-{\tilde u}=g_1,  \hfill & \cr
i\lambda {\tilde u}-(x^{\gamma}u_x)_x+\alpha v=g_2, \hfill & \cr
i\lambda v-{\tilde v}=g_3\hfill & \cr
i\lambda {\tilde v}-(x^{\gamma}v_x)_x+\alpha u=g_4, \hfill & \cr
i\lambda \varphi+(\varsigma^{2}+\omega)\varphi-{\tilde u}(0)\vartheta(\varsigma)=g_5.\hfill & \cr}\right.
\label{e35e}
\end{equation}
Inserting $(\ref{e35e})_1$ and $(\ref{e35e})_3$ into $(\ref{e35e})_2$ and $(\ref{e35e})_4$, we get
\begin{equation}
\left\{\matrix{-\lambda^2 u-(x^{\gamma}u_x)_x+\alpha v=(g_2+i\lambda g_1), \hfill & \cr
-\lambda^2 v-(x^{\gamma}v_x)_x+\alpha u=(g_4+i\lambda g_3), \hfill & \cr}\right.
\label{mk}
\end{equation}
Solving system (\ref{mk}) is equivalent to finding $u\in H_{\gamma}^{2}\cap H_{0, \gamma}^{1}(0,1)$ and $v\in H_{\gamma}^{2}\cap H_{00, \gamma}^{1}(0,1)$
such that
\begin{equation}
\left\{\matrix{\Int_{0}^{1}(-\lambda^{2}u {\overline w}- (x^{\gamma} u_x)_x {\overline w}+\alpha v{\overline w})\, dx=
\Int_{0}^{1}(g_2+i\lambda g_1){\overline w}\, dx,\hfill & \cr
\Int_{0}^{1}(-\lambda^{2}u {\overline y}- (x^{\gamma} u_x)_x {\overline y}+\alpha u{\overline y})\, dx=
\Int_{0}^{1}(g_4+i\lambda g_3){\overline y}\, dx,\hfill & \cr}\right.
\label{mk1}
\end{equation}
for all $w\in H_{0, \gamma}^{1}(0,1)$ and $y\in H_{00, \gamma}^{1}(0,1)$. Then, we get
\begin{equation}
\left\{\matrix{
\Int_{0}^{1}(-\lambda^{2} u {\overline w}+ x^{\gamma} u_x  {\overline w_x}+\alpha v{\overline w})\, dx
+i\lambda {\tilde\zeta}  u(0)\ {\overline w(0)}\, \hfill &\cr
=\Int_{0}^{1}(g_2+i\lambda g_1){\overline w}\, dx
-\zeta \Int_{-\infty}^{+\infty} \Frac{\vartheta(\varsigma)}{\varsigma^{2}+\omega+i\lambda}f_5(\varsigma)  \, d\varsigma{\overline w}(0)
+{\tilde\zeta}  g_1(1) {\overline w(0)},\hfill & \cr
\Int_{0}^{1}(-\lambda^{2}v {\overline y}+ x^{\gamma} v_x {\overline y_x}+\alpha u{\overline y})\, dx=
\Int_{0}^{1}(g_4+i\lambda g_3){\overline y}\, dx.\hfill & \cr}\right.
\label{mk2}
\end{equation}
We can rewrite (\ref{mk2}) as
\begin{equation}
{\cal B}((u, v), (w, y))=l(w, y),\quad \forall (w, y)\in H_{0, \gamma}^{1}(0,1)\times H_{00, \gamma}^{1}(0,1),
\label{e2xx5}
\end{equation}
where
$$
{\cal B}(u, v), (w, y)={\cal B}_1(u, v), (w, y)+{\cal B}_2(u, v), (w, y)
$$
with
$$
\left\{\matrix{{\cal B}_1(u, v), (w, y)=\Int_{0}^{1}  (x^{\gamma}(u_x {\overline w_x}+v_x {\overline y_x})
+\alpha (v{\overline w}+ u{\overline y}))\, dx+i\lambda {\tilde\zeta} u(0)\ {\overline w(0)},\hfill \cr
{\cal B}_2(u, v), (w, y)=-\Int_{0}^{1}\lambda^{2}(u {\overline w}+v {\overline y})\, dx,\hfill \cr}\right.
\leqno{(*)}
$$
and
$$
\matrix{l(w, y)=\Int_{0}^{1}(g_2+i\lambda g_1){\overline w}\, dx
-\zeta\Int_{-\infty}^{+\infty} \Frac{\vartheta(\varsigma)}{\varsigma^{2}+\omega+i\lambda}g_3(\varsigma)  \, d\varsigma\ {\overline w}(0)\hfill &\cr
+\rho(i\lambda+\omega)^{\tau-1}g_1(0){\overline w}(0)+\Int_{0}^{1}(g_4+i\lambda g_3){\overline y}\, dx.
&\cr}
$$
Let $(H_{0, \gamma}^{1}(0,1)\times H_{00, \gamma}^{1}(0,1))'$ be the dual space of $H_{0, \gamma}^{1}(0,1)\times H_{00, \gamma}^{1}(0,1)$.
Let us define the following operators
$$
\matrix{B:H_{0, \gamma}^{1}\times H_{00, \gamma}^{1}\rightarrow (H_{0, \gamma}^{1}\times H_{00, \gamma}^{1})'\hfill &\cr
\quad (u, v)\mapsto B (u, u) &\cr}
\matrix{B_i:H_{0, \gamma}^{1}\times H_{00, \gamma}^{1}(0,1)\rightarrow (H_{0, \gamma}^{1}\times H_{00, \gamma}^{1}(0,1))'\quad i\in \{1, 2\}\hfill &\cr
\quad (u, v)\mapsto B_i (u, v) &\cr}
\leqno{(**)}
$$
such that
$$
\matrix{(B (u, v)) (w, y)={\cal B}((u, v), (w, y)), \ \forall (w, y)\in H_{0, \gamma}^{1}\times H_{00, \gamma}^{1}(0,1),\hfill &\cr
(B_i u) w={\cal B}_i(u, w), \ \forall (w, y)\in H_{0, \gamma}^{1}\times H_{00, \gamma}^{1}(0,1), i\in \{1, 2\}.\hfill &\cr}
\leqno{(***)}
$$
We need to prove that the operator $B$ is an isomorphism. For this aim, we divide the proof into three steps:

\noindent
{\bf Step 1.}
In this step, we want to prove that the operator $B_1$ is an isomorphism. For this aim, it is easy to see that
${\cal B}_1$ is sesquilinear, continuous form on $H_{0, \gamma}^{1}(0,1)\times H_{00, \gamma}^{1}(0,1)$.
Furthermore
$$
\matrix{\Re {\cal B}_1((u, v), (u, v))&=&\| x^{\gamma/2}u_x\|_{2}^2+\|x^{\gamma/2}v_x\|_{2}^2
+\alpha\Int_{0}^{1}(u{\overline v}+v{\overline u})\, dx\hfill \cr
&&+\rho\lambda\Re \left(i(i\lambda+\omega)^{\tau-1}\right)|u(1)|^2\hfill \cr
&\geq & c(\|x^{\gamma/2}u_x\|_{2}^2+\|x^{\gamma/2}v_x\|_{2}^2),\hfill \cr}
$$
where we have used the fact that
$$
\rho\lambda\Re \left(i(i\lambda+\omega)^{\tau-1}\right)=\zeta\lambda^2\Int_{-\infty}^{+\infty}\Frac{\vartheta(\varsigma)^2}{\lambda^2+(\omega+\varsigma^2)^2}\, d\varsigma> 0.
$$
Thus ${\cal B}_1$ is coercive. Then, from $(**)$ and Lax-Milgram theorem, the operator $B_1$ is
an isomorphism.

\noindent
{\bf Step 2.}
In this step, we want to prove that the operator $B_2$ is compact. For this aim, from $(*)$ and $(***)$, we have
$$
|{\cal B}_2((u, v), (w, y))|\leq c\|(u, v)\|_{L^{2}(0, 1)}\|(w, y)\|_{L^{2}(0, 1)},
$$
and consequently, using the compact embedding from $H_{0, \gamma}^{1}(0,1)\times H_{00, \gamma}^{1}(0,1))$
to $L^{2}(0, 1)\times L^{2}(0, 1)$
we deduce that $B_2$ is a compact operator.
Therefore, from the above steps, we obtain that the operator $B =B_1+ B_2$ is a Fredholm operator of index zero. Now,
following Fredholm alternative, we still need to prove that the operator $B$ is injective to obtain that the operator $B$
is an isomorphism.

\noindent
{\bf Step 3.}
Let $(u, v)\in ker(B)$, then
\begin{equation}
{\cal B}(u, v), (w, y))=0\quad \forall (w, y)\in H_{0, \gamma}^{1}(0,1)\times H_{00, \gamma}^{1}(0,1).
\label{eaz}
\end{equation}
In particular for $(w, y)=(u, v)$, it follows that
$$
\matrix{\lambda^{2}(\|u\|_{L^2(0, 1)}^2+\|v\|_{L^2(0, 1)}^2)
-i\rho\lambda(i\lambda+\omega)^{\tau-1}|u(0)|^2=\hfill & \cr
\qquad\|x^{\gamma/2}u_x\|_{L^2(0, 1)}^{2}+\|x^{\gamma/2}v_x\|_{L^2(0, 1)}^{2}
+\alpha\Int_{0}^{1}(v{\overline u}+v{\overline u})\, dx. & \cr}
$$
Hence, we have
\begin{equation}
u(0)=0.
\label{e22xx1}
\end{equation}
From (\ref{eaz}), we obtain
\begin{equation}
(x^{\gamma}u_x)(0)=0
\label{e23xx1}
\end{equation}
and then
\begin{equation}
\left\{\matrix{
\lambda^2 u+(x^{\gamma}u_x)_x-\alpha v=0, \hfill & \cr
\lambda^2 v+(x^{\gamma}v_x)_x-\alpha u=0, \hfill & \cr
u(1)=v(0)=v(1)=0, \hfill & \cr
u(0)=(x^{\gamma}u_x)(0)=0 \hfill &,\hfill \cr
}\right.
\label{aeq44b}
\end{equation}
Then, according to Lemma \ref{lem10}, we deduce that $(u, v)=(0, 0)$ and consequently $Ker(B)=\{0\}$.
Finally, from Step 3 and Fredholm alternative, we deduce that the operator $B$ is isomorphism. It is easy
to see that the operator $l$ is a antilinear and continuous form on $H_{0, \gamma}^{1}\times H_{00, \gamma}^{1}(0,1)$. Consequently,
(\ref{e2xx5}) admits a unique solution $(u, v)\in H_{0, \gamma}^{1}(0,1)\times H_{00, \gamma}^{1}(0,1)$.
By using the classical elliptic regularity, we deduce that $U\in D({\cal P})$ is a unique solution
of (\ref{e35eee}). Hence $i\lambda- {\cal P}$ is surjective for all $\lambda\in \R^*$.\\

\noindent
{\bf Case 2: $\lambda=0$ and $\omega\not=0$.} Using Lax-Milgram Lemma, we obtain the result.

Taking account of Lemmas \ref{lem10}, \ref{lem11}
and from Theorem \ref{thm14} the $C_0$-semigroup $e^{t{\cal P}}$ is strongly stable in ${\cal H}$.

\hfill$\diamondsuit$\\

\section{Lack of exponential stability}
This section will be devoted to the study of the lack of exponential decay of solutions associated with the system
$(P')$.
\begin{proposition}
The $C_0$-semigroup of contractions $S(t) = e^{{\cal P} t}$ associated with (\ref{e14}) is not exponentially stable.
\label{pro1}
\end{proposition}
{\bf Proof.}
Let $\eta_n$ be an eigenvalue of ${\cal A}u=-(x^{\gamma}u_x)_x$ in $H_{0, \gamma}^{1}(0, 1)$ corresponding to the normalized
eigenfunction $e_n$, and
$$
U_n=\Frac{1}{\sqrt{2}}\left(0, 0, \frac{e_n}{i\sqrt{\eta_n}}, e_n, 0\right)^T.
$$
Then a straightforward computation gives
$$
\|U_n\|_{{\cal H}}=1, \qquad \|(i\sqrt{\eta_n}-{\cal P})U_n\|_{{\cal H}}^2=\Frac{\alpha^2}{2\eta_n}\rightarrow 0.
$$
This shows that the resolvent of ${\cal P}$ is not uniformly bounded on the imaginary axis. Following
{\bf\cite{12}} and {\bf\cite{13}}, the system $(P')$ is not uniformly and exponentially stable in the energy
space ${\cal H}$.\\

\noindent
\subsection{Precise spectral analysis}
We aim to show that an
infinite number of eigenvalues of ${\cal P}$ approach the imaginary
axis which prevents the system (\ref{eqq100}) from being exponentially
stable. Indeed we first compute the characteristic equation that
gives the eigenvalues of  ${\cal P}$.
Let $\lambda$ be an eigenvalue of ${\cal P}$ with associated
eigenvector $U=(u,{\tilde u}, v, {\tilde v}, \varphi)^{T}$. Then ${\cal P}U = \lambda U$ is
equivalent to
\begin{equation}
\left\{\matrix{\lambda u-{\tilde u}=0,  \hfill & \cr
\lambda {\tilde u}-(x^{\gamma}u_x)_x+\alpha v=0, \hfill & \cr
\lambda v-{\tilde v}=0\hfill & \cr
\lambda {\tilde u}-(x^{\gamma}v_x)_x+\alpha u=0, \hfill & \cr
(\lambda+\varsigma^2+\omega)\varphi-{\tilde u}(0)\vartheta(\varsigma)=0,
}\right.
\label{066}
\end{equation}
with boundary conditions
\begin{equation}
\left\{\matrix{u(1)=v(0)=v(1)=0, \hfill & \cr
 (x^{\gamma}u_x)(0)=\rho \lambda(\lambda+\omega)^{\tau-1}u(0).\hfill &\cr}\right.
\label{e556aa}
\end{equation}
Inserting $(\ref{066})_1$ into $(\ref{066})_2$ and $(\ref{066})_3$ into $(\ref{066})_4$, we get
\begin{equation}
\left\{\matrix{
\lambda^2 u-(x^{\gamma}u_x)_x+\alpha v=0, \hfill & \cr
\lambda^2 v-(x^{\gamma}v_x)_x+\alpha u=0, \hfill & \cr
u(1)=v(0)=v(1)=0, \hfill & \cr
 (x^{\gamma}u_x)(0)=\rho \lambda(\lambda+\omega)^{\tau-1}u(0).
}\right.
\label{e331cc}
\end{equation}
Let us set
\begin{equation}
\left\{\matrix{ \phi=u+v, \hfill \cr
 {\psi}=u-v.
\hfill &\cr}\right.
\label{e771}
\end{equation}
Then, we obtain
 \begin{equation}
 \left\{\matrix{(\lambda^2+\alpha)\phi-(x^{\gamma}\phi_x)_x=0,\hfill &\cr
(\lambda^2-\alpha)\psi-(x^{\gamma}\psi_x)_x=0.
\hfill &\cr}\right.
\label{e775}
 \end{equation}
The solution of equations (\ref{e775}) is given by
\begin{equation}
\left\{\matrix{\phi(x)=c_{1}\Phi_{+}+c_{2}\Phi_{-},\hfill &\cr
\psi(x)={\tilde c}_{1}\Phi_{++}+{\tilde c}_{2}\Phi_{--},\hfill &\cr}\right.
\end{equation}
where $\Phi_{+}, \Phi_{-}, \Phi_{++}$ and $\Phi_{--}$ are defined by
$$
\left\{\matrix{
\Phi_{+}(x):=x^{\frac{1-\gamma}{2}}J_{\nu_\gamma}\left(\frac{2}{2-\gamma}i{\tilde\lambda}
x^{\frac{2-\gamma}{2}}\right),\ \
\Phi_{-}(x):=x^{\frac{1-\gamma}{2}}J_{-\nu_\gamma}\left(\frac{2}{2-\gamma}i{\tilde\lambda}
x^{\frac{2-\gamma}{2}}\right),\hfill &\cr\cr
\Phi_{++}(x):=x^{\frac{1-\gamma}{2}}J_{\nu_\gamma}\left(\frac{2}{2-\gamma}i{\tilde{\tilde\lambda}}
x^{\frac{2-\gamma}{2}}\right),\ \
\Phi_{--}(x):=x^{\frac{1-\gamma}{2}}J_{-\nu_\gamma}\left(\frac{2}{2-\gamma}i{\tilde{\tilde\lambda}}
x^{\frac{2-\gamma}{2}}\right),\hfill &\cr}\right.
$$
where
$$
{\tilde\lambda}=\sqrt{\lambda^2+\alpha},\ \ {\tilde{\tilde\lambda}}=\sqrt{\lambda^2-\alpha},\ \
\nu_\gamma=\Frac{1-\gamma}{2-\gamma}.
$$
Then
$$
\left\{\matrix{u(x)=\Frac{1}{2}(c_{1}\Phi_{+}+c_{2}\Phi_{-}+{\tilde c}_{1}\Phi_{++}+{\tilde c}_{2}\Phi_{--}),
\hfill &\cr
v(x)=\Frac{1}{2}(c_{1}\Phi_{+}+c_{2}\Phi_{-}-{\tilde c}_{1}\Phi_{++}-{\tilde c}_{2}\Phi_{--}).
\hfill &\cr
}\right.
$$
By substituting the series expansions of $J_{\nu_{\gamma}}$ and
$J_{-\nu_{\gamma}}$, we obtain
$$
\matrix{\Phi_{+}(x)=\Sum_{m=0}^{\infty}{\tilde c}_{\nu_\gamma, m}^{+}x^{1-\gamma+(2-\gamma)m},\quad
\Phi_{-}(x)=\Sum_{m=0}^{\infty}{\tilde c}_{\nu_\gamma, m}^{-}x^{(2-\gamma)m}\hfill &\cr
\Phi_{++}(x)=\Sum_{m=0}^{\infty}{\tilde e}_{\nu_\gamma, m}^{+}x^{1-\gamma+(2-\gamma)m},\quad
\Phi_{--}(x)=\Sum_{m=0}^{\infty}{\tilde e}_{\nu_\gamma, m}^{-}x^{(2-\gamma)m}\hfill &\cr}
$$
with
$$
\left\{\matrix{{\tilde c}_{\nu_\gamma, m}^{+}=c_{\nu_\gamma, m}^{+}\left(\Frac{2}{2-\gamma}i{\tilde\lambda}\right)^{2m+\nu_\gamma},\quad
{\tilde c}_{\nu_\gamma, m}^{-}=c_{\nu_\gamma, m}^{-}\left(\Frac{2}{2-\gamma}i{\tilde\lambda}\right)^{2m-\nu_\gamma}\hfill & \cr
{\tilde e}_{\nu_\gamma, m}^{+}=c_{\nu_\gamma, m}^{+}\left(\Frac{2}{2-\gamma}i{\tilde{\tilde\lambda}}\right)^{2m+\nu_\gamma},\quad
{\tilde e}_{\nu_\gamma, m}^{-}=c_{\nu_\gamma, m}^{-}\left(\Frac{2}{2-\gamma}i{\tilde{\tilde\lambda}}\right)^{2m-\nu_\gamma}.\hfill & \cr}\right.
$$
One can readily verify that $\Phi_{+}, \Phi_{++},\Phi_{-}, \Phi_{--}$ belongs to $\in H_{\gamma}^1(0, 1)$. indeed,
$$
\matrix{{\Phi}_{+}(x)\sim_0 {\tilde c}_{\nu_\gamma, 0}^{+}x^{1-\gamma},\quad
x^{\gamma/2}{\Phi}'_{+}(x)\sim_0 (1-\gamma){\tilde c}_{\nu_\gamma, 0}^{+}x^{-\gamma/2},\hfill \cr
{\Phi}_{-}(x)\sim_0 {\tilde c}_{\nu_\gamma, 0}^{-},\quad
x^{\gamma/2}{\Phi}'_{-}(x)\sim_0 (2-\gamma){\tilde c}_{\nu_\gamma, 0}^{-}x^{1-\gamma/2},\hfill \cr}
$$
where we have used the following relation
\begin{equation}
x J'_{\nu}(x)=\nu J_{\nu}(x)-x J_{\nu+1}(x).
\label{e7112}
\end{equation}
With the boundary conditions
$$
\left\{\matrix{u(1)=v(0)=v(1)=0, \hfill & \cr
(x^{\gamma}u_{x})(0)=\rho \lambda(\lambda+\omega)^{\tau-1}u(0),\hfill & \cr}\right.
$$
we have
\begin{equation}
M\pmatrix{c_{1}\cr c_{2}\cr {\tilde c}_1\cr {\tilde c}_2\cr
}=\pmatrix{0\cr 0\cr 0\cr 0\cr},
\label{mat1}
\end{equation}
where
$$
M=\pmatrix{(1-\gamma){\tilde c}_{\nu_\gamma, 0}^{+} &(1-\gamma){\tilde{\tilde c}}_{\nu_\gamma, 0}^{+}&
-\rho \lambda(\lambda+\omega)^{\tau-1}{\tilde c}_{\nu_\gamma, 0}^{-}&-\rho \lambda(\lambda+\omega)^{\tau-1}{\tilde{\tilde c}}_{\nu_\gamma, 0}^{-}\cr
0 & 0&{\tilde c}_{\nu_\gamma, 0}^{-}& -{\tilde{\tilde c}}_{\nu_\gamma, 0}^{-} \cr
J_{\nu_{\gamma}}(K_{+})&J_{-\nu_{\gamma}}(K_{+})&J_{\nu_{\gamma}}(K_{-})&J_{-\nu_{\gamma}}(K_{-})\cr
J_{\nu_{\gamma}}(K_{+})&J_{-\nu_{\gamma}}(K_{+})&-J_{\nu_{\gamma}}(K_{-})&-J_{-\nu_{\gamma}}(K_{-})\cr}.
$$
System (\ref{e331cc}) admits a non trivial solution if and only if $det(M)=0$, which is equivalent to requiring that
the eigenvalues of ${\cal P}$ be zeros of the function $f$ defined by
\begin{equation}
\matrix{f(\lambda)&=&
C_{\gamma}(K_+)J_{-\nu_{\gamma}}(K_{+})J_{\nu_{\gamma}}(K_{-})+C_{\gamma}(K_-)J_{-\nu_{\gamma}}(K_{-})J_{\nu_{\gamma}}(K_{+})\hfill \cr
&& +2\rho \lambda(\lambda+\omega)^{\tau-1}J_{\nu_{\gamma}}(K_{+})J_{\nu_{\gamma}}(K_{-}),\hfill \cr}
\label{sszz}
\end{equation}
where
$$
K{\pm}=\frac{2}{2-\gamma}i\sqrt{\lambda^2\pm\alpha},\ \ C_{\gamma}(K)=(1-\gamma)\Frac{\Gamma(1-\nu_{\gamma})}{\Gamma(1+\nu_{\gamma})}\left(\Frac{K}{2}\right)^{2\nu_{\gamma}}.
$$
Our aim is to show, by means of Rouch\'e's Theorem, that there exists a subsequence of eigenvalues whose real parts converge to zero.

In what follows, since ${\cal P}$ is dissipative, we investigate the asymptotic behavior of the eigenvalues $\lambda$ of ${\cal P}$ with sufficiently
large modulus in the vertical strip $\Sigma=\{\lambda\in \C: -\alpha_0\leq \Re(\lambda)\leq 0\}$,
where $\alpha_0>0$ is sufficiently large. For such eigenvalues, we observe that $\Phi_{+}, \Phi_{-}, \Phi_{++}, \Phi_{--}$ remain bounded.
Our main result in this section is the following.
\begin{lemma}
The large eigenvalues of the dissipative operator ${\cal P}$ are simple and can
be split into two families $(\lambda_{j, k})_{k\in \Z, |k|\geq N}, j = 1, 2, (N\in  \N$, chosen large enough).
Moreover, the following asymptotic expansions for the eigenvalues hold:

\noindent
$\bullet$ If $2\nu_{\gamma}>\tau$, then
$$
\lambda_{1, k}=\frac{1}{iq}\left(k\pi-\nu_{\gamma}\frac{\pi}{2}+\frac{3}{4}\pi\right)-\Frac{\tilde C}{q^{\tau+1}}\Frac{\sin(\nu_{\gamma}\pi)}{(k\pi)^{2\nu_{\gamma}-\tau}}i^{1-\tau}+O\left(\Frac{1}{k}\right), \quad k\in\Z.
$$
$$
\lambda_{2, k}=\frac{1}{iq}\left(k\pi+\nu_{\gamma}\frac{\pi}{2}+\frac{3}{4}\pi\right)+i\Frac{a_1}{q k\pi}
-\Frac{{\tilde C}}{8}\Frac{\alpha^2 q^{3-\tau}}{\sin(\nu_{\gamma}\pi)}\Frac{i^{1-\tau}}{(k\pi)^{2+2\nu_{\gamma}-\tau}}+O\left(\Frac{1}{k^{3}}\right).
$$
$$
|k|^{2\nu_{\gamma}-\tau}\Re\lambda_{1, k}=\upsilon=-\Frac{\tilde C}{q^{\tau+1}}\Frac{\sin(\nu_{\gamma}\pi)\sin\frac{\pi}{2}\tau}{\pi^{2\nu_{\gamma}-\tau}}.
$$
$$
|k|^{2+2\nu_{\gamma}-\tau}\Re\lambda_{2, k}=\upsilon=-\Frac{{\tilde C}}{8}\Frac{\alpha^2 q^{3-\tau}}{\pi^{2+2\nu_{\gamma}-\tau}}\Frac{\sin\frac{\tau}{2}\pi}{\sin(\nu_{\gamma}\pi)}
$$

\noindent
$\bullet$ If $2\nu_{\gamma}<\tau$, then
$$
\lambda_{1, k}=\frac{1}{iq}\left(k\pi+\nu_{\gamma}\frac{\pi}{2}+\frac{3}{4}\pi\right)
-\frac{2 q^{\tau-1}\sin\nu_\gamma\pi}{{\tilde C}(k\pi)^{\tau-2\nu_\gamma}}i^{\tau-1}+o\bigl(k^{-(\tau-2\nu_{\gamma})}\bigr),
$$
where
$$
|k|^{\tau-2\nu_{\gamma}}\Re\lambda_{1, k}=\upsilon=-\frac{2 q^{\tau-1}\sin\nu_\gamma\pi\sin\frac{\tau}{2}\pi}{{\tilde C}\pi^{\tau-2\nu_\gamma}}.
$$
$$
\lambda_{2, k}=\frac{1}{iq}\left(k\pi+\nu_{\gamma}\frac{\pi}{2}+\frac{3}{4}\pi\right)
-\frac{a_1}{l\pi k}-\frac{{\tilde C}q^{3-\tau} \alpha^2}{2\sin\nu_\gamma\pi\,(k\pi)^{2+2\nu_{\gamma}-\tau}}i^{1-\tau}+o\bigl(k^{-(2+2\nu_{\gamma}-\tau)}\bigr)
$$
where
$$
|k|^{2+2\nu_{\gamma}-\tau}\Re\lambda_{2, k}=\upsilon=-\frac{{\tilde C} q^{3-\tau}\alpha^2\sin\frac{\tau}{2}\pi}{2\sin\nu_\gamma\pi\,\pi^{2+2\nu_{\gamma}-\tau}}.
$$

\noindent
$\bullet$ If $2\nu_{\gamma}=\tau$, then
$$
\matrix{\lambda_{1, k}&=&\left[-\frac{1}{q}(k\pi+\nu_{\gamma}\frac{\pi}{2}+\frac{3}{4}\pi)+\frac{a_1}{q\pi}\frac1k+\left(-\frac{a_1}{q\pi^2}\left(\frac{3\pi}{4}+\frac{\pi\tau}{4}\right)
+\frac{q b_1^2}{2\pi^2}\frac{\cos\frac{\pi\tau}{2}}{\sin\frac{\pi\tau}{2}}
\left(2+\frac{{\tilde C}}{q^\tau}\right)\right)\frac{1}{k^2}\right]i\hfill
\cr &&+\Frac{q^{1-\tau} b_1^2 {\tilde C}}{2\pi^2}\frac{1}{k^2}+O(k^{-3}),\hfill  \cr}
$$
where
$$
|k|^{2}\Re\lambda_{1, k}\sim \upsilon_1=-\Frac{q^{3-\tau}\alpha^2 {\tilde C}}{8\pi^2}.
$$
$$
\lambda_{2, k}=-i\frac{2-\gamma}{4}(2k\pi-\Frac{\pi}{2}+\theta)-\frac{2-\gamma}{4}\ln\Frac{2+r}{\sqrt{4+r^2+4r\cos
2\nu_\gamma\pi}}+O(k^{-1}),\quad k\in\Z.
$$
where $\theta$ is such that
$$
\left\{\matrix{\cos\theta&=&\Frac{(2+r)\cos\nu_\gamma\pi}{\sqrt{4+r^2+4r\cos2\nu_\gamma\pi}},\hfill\cr
\sin\theta&=&\Frac{(2-r)\sin\nu_\gamma\pi}{\sqrt{4+r^2+2r\cos2\nu_\gamma\pi}}.\hfill\cr}\right.
$$
$$
\lambda_{j, k}=\overline{\lambda_{j, -k}} \hbox{ if } k\leq -N,
$$
Moreover for all $|k|\geq N$, the eigenvalues $\lambda_{j, k}$ are simple.
\label{ll12}
\end{lemma}

\noindent
{\bf Proof.}

\noindent
{\bf Step 1.}
We will use the following classical development (see
{\bf\cite{lebed}} p. 122, (5.11.6)): for all $\delta>0$, the
following development holds when $|\arg z|< \pi-\delta$:
\begin{equation}
\matrix{J_{\nu}(z)&=&\left(\Frac{2}{\pi z}\right)^{1/2}\left[\cos\chi_{\nu}
-\Frac{(4\nu^2-1)}{8 z}\sin\chi_{\nu}
-\Frac{(4\nu^2-1)(4\nu^2-9)}{2(8 z)^2}\cos\chi_{\nu}
+O\left(\frac{1}{|z|^3}\right)\right] \hfill \cr
&=&\left(\Frac{2}{\pi z}\right)^{1/2}\left[\cos\chi_{\nu}-\Frac{a_1}{z}\sin\chi_{\nu}-\Frac{a_2}{z^2}\cos\chi_{\nu}
+O\left(\frac{1}{|z|^3}\right)\right]\hfill \cr}
\label{029}
\end{equation}
where $\chi_{\nu}=\left(z-\nu\frac{\pi}{2}-\frac{\pi}{4}\right)$.
Moreover, for $\lambda$ large enough, we have
\begin{equation}
\matrix{{\tilde\lambda}=\sqrt{\lambda^2+\alpha}&=&\lambda+\Frac{\alpha}{2\lambda}-\Frac{\alpha^2}{8\lambda^3}+\Frac{O(\alpha^3)}{\lambda^5},\hfill \cr
{\tilde{\tilde\lambda}}=\sqrt{\lambda^2-\alpha}&=&\lambda-\Frac{\alpha}{2\lambda}-\Frac{\alpha^2}{8\lambda^3}+\Frac{O(\alpha^3)}{\lambda^5},\hfill \cr
{\tilde\lambda}+{\tilde{\tilde\lambda}}&=&2\lambda-\Frac{\alpha^2}{4\lambda^3}+\Frac{O(\alpha^2)}{\lambda^4},\hfill \cr
{\tilde\lambda}-{\tilde{\tilde\lambda}}&=&\Frac{\alpha}{\lambda}+\Frac{O(\alpha^2)}{\lambda^4},\hfill \cr
\lambda(\lambda+\omega)^{\tau-1}&=&\lambda^{\tau}+\Frac{(\tau-1)\omega}{\lambda^{1-\tau}}
+\Frac{(\tau-1)(\tau-2)\omega^2}{2\lambda^{2-\tau}}+O\left(\Frac{1}{\lambda^{3-\tau}}\right).\hfill \cr
}
\label{sq2}
\end{equation}
Set
\begin{equation}
\matrix{
I_{1, \nu}&=&\cos\chi_{-\nu,+}-\Frac{a_1}{l\lambda}\sin\chi_{-\nu,+}
-\Frac{a_2}{l^2\lambda^2}\cos\chi_{-\nu,+}+O\left(\Frac{1}{\lambda^3}\right),\hfill\cr
I_{2, \nu}&=&\cos\chi_{\nu,-}-\Frac{a_1}{l\lambda}\sin\chi_{\nu,-}
-\Frac{a_2}{l^2\lambda^2}\cos\chi_{\nu,-}+O\left(\Frac{1}{\lambda^3}\right),\hfill\cr
I_{3, \nu}&=&\cos\chi_{\nu,+}-\Frac{a_1}{l\lambda}\sin\chi_{\nu,+}
-\Frac{a_2}{l^2\lambda^2}\cos\chi_{\nu,+}+O\left(\Frac{1}{\lambda^3}\right),\hfill\cr
I_{4, \nu}&=&\cos\chi_{-\nu,-}-\Frac{a_1}{l\lambda}\sin\chi_{-\nu,-}
-\Frac{a_2}{l^2\lambda^2}\cos\chi_{-\nu,-}+O\left(\Frac{1}{\lambda^3}\right),\hfill \cr}
\label{eqqq1}
\end{equation}
where
\begin{equation}
\matrix{
\cos\chi_{-\nu,+}
&=\cos A^*-\Frac{b_1}{\lambda}\sin A^*
+\Frac{b_2}{\lambda^2}\cos A^*+O\left(\Frac{1}{\lambda^3}\right),\hfill \cr
\sin\chi_{-\nu,+}
&=\sin A^*+\Frac{b_1}{\lambda}\cos A^*
+\Frac{b_2}{\lambda^2}\sin A^*+O\left(\Frac{1}{\lambda^3}\right),\hfill \cr
\cos\chi_{\nu,-}
&=\cos A^{**}+\Frac{b_1}{\lambda}\sin A^{**}
+\Frac{b_2}{\lambda^2}\cos A^{**}+O\left(\Frac{1}{\lambda^3}\right),\hfill \cr
\sin\chi_{\nu,-}
&=\sin A^{**}-\Frac{b_1}{\lambda}\cos A^{**}
+\Frac{b_2}{\lambda^2}\sin A^{**}+O\left(\Frac{1}{\lambda^3}\right),\hfill \cr
\cos\chi_{-\nu,-}
&=\cos A^*+\Frac{b_1}{\lambda}\sin A^*
+\Frac{b_2}{\lambda^2}\cos A^*+O\left(\Frac{1}{\lambda^3}\right),\hfill \cr
\sin\chi_{-\nu,-}
&=\sin A^*-\Frac{b_1}{\lambda}\cos A^*
+\Frac{b_2}{\lambda^2}\sin A^*+O\left(\Frac{1}{\lambda^3}\right),\hfill \cr
\cos\chi_{\nu,+}
&=\cos A^{**}-\Frac{b_1}{\lambda}\sin A^{**}
+\Frac{b_2}{\lambda^2}\cos A^{**}+O\left(\Frac{1}{\lambda^3}\right),\hfill \cr
\sin\chi_{\nu,+}
&=\sin A^{**}+\Frac{b_1}{\lambda}\cos A^{**}
+\Frac{b_2}{\lambda^2}\sin A^{**}+O\left(\Frac{1}{\lambda^3}\right),\hfill \cr}
\label{eqqq2}
\end{equation}
where
$$
\matrix{A^*=\frac{2}{2-\gamma}i\lambda+\nu_{\gamma}\frac{\pi}{2}-\frac{\pi}{4},\hfill &\cr
A^{**}=\frac{2}{2-\gamma}i\lambda-\nu_{\gamma}\frac{\pi}{2}-\frac{\pi}{4}.\hfill &\cr}
$$
where
$$
\left\{\matrix{a_1=\Frac{(4\nu_{\gamma}^2-1)}{8},\
a_2=\Frac{(4\nu_{\gamma}^2-1)(4\nu_{\gamma}^2-9)}{128},\ \hfill & \cr
b_1=\Frac{iq\alpha}{2},\ b_2=\Frac{q^2\alpha^2}{8},\ q=\frac{2}{2-\gamma}, \ l=iq,\hfill &\cr}\right.
$$
Inserting (\ref{eqqq1}) in (\ref{sszz}) and using (\ref{eqqq2}) we get
$$
\matrix{f(\lambda)&=&(1-\gamma)\Frac{\Gamma(1-\nu_{\gamma})}{\Gamma(1+\nu_{\gamma})}\Frac{(iq)^{2\nu_{\gamma}}}{2^{2\nu_{\gamma}}}\left(\sqrt{\lambda^2+\alpha}^{2\nu_{\gamma}} J_{-\nu_{\gamma}}(K_{+})J_{\nu_{\gamma}}(K_{-})+\sqrt{\lambda^2-\alpha}^{2\nu_{\gamma}}J_{-\nu_{\gamma}}(K_{-})J_{\nu_{\gamma}}(K_{+})\right.\hfill \cr
&& \left.+\Frac{2^{1+2\nu_{\gamma}}\rho\Gamma(1+\nu_{\gamma})}{(1-\gamma)\Gamma(1-\nu_{\gamma})(iq)^{2\nu_{\gamma}}} \lambda(\lambda+\omega)^{\tau-1}J_{\nu_{\gamma}}(K_{+})J_{\nu_{\gamma}}(K_{-})\right)\hfill \cr
&=&(1-\gamma)\Frac{\Gamma(1-\nu_{\gamma})}{\Gamma(1+\nu_{\gamma})}\Frac{(iq)^{2\nu_{\gamma}}}{2^{2\nu_{\gamma}}}\left(\Frac{2}{\pi K_+}\right)^{1/2}\left(\Frac{2}{\pi K_-}\right)^{1/2}{\tilde f}(\lambda),\hfill \cr}
$$
where
\begin{equation}
\matrix{{\tilde f}(\lambda)&=&
2\lambda^{2\nu_{\gamma}}\cos A^*\cos A^{**}-\Frac{2a_1}{l\lambda^{1-2\nu_{\gamma}}}
\sin(A^*+A^{**})\hfill  \cr
&&+\Frac{1}{\lambda^{2-2\nu_{\gamma}}}\left(4\left(b_2-\Frac{a_2}{l^2}\right)
\cos A^*\cos A^{**}+ 2\left(\Frac{a_1^2}{l^2}-b_1^2\right)\sin A^*\sin A^{**}\right)\hfill \cr
&&+\Frac{2\nu_{\gamma}\alpha b_1}{\lambda^{3-2\nu_{\gamma}}}\sin(A^{**}-A^{*})+O\left(\Frac{1}{\lambda^{4-2\nu_{\gamma}}}\right)\hfill \cr
&&+C\left[\lambda^\tau\cos^2 A^{**}
+ \lambda^{\tau-1}
\left(-\Frac{a_1}{l}\sin(2A^{**})
+(\tau-1)\omega\cos^2 A^{**}\right)\right.\hfill \cr
&&+\lambda^{\tau-2}\Bigg(2\left(b_2-\Frac{a_2}{l^2}\right)
\cos^2 A^{**}+\left(\Frac{a_1^2}{l^2}-b_1^2\right)
\sin^2 A^{**}\hfill \cr
&&\left.-(\tau-1)\omega\Frac{a_1}{l}\sin(2A^{**})+\Frac{(\tau-1)(\tau-2)}{2}\omega^2\cos^2 A^{**}\Bigg)\right]\hfill \cr
&&+O(\lambda^{\tau-3}).\hfill \cr
}
\label{030}
\end{equation}
where
\begin{equation}
C=\Frac{2^{2\nu_{\gamma}+1}\rho}{(1-\gamma)(iq)^{2\nu_{\gamma}}}\Frac{\Gamma(1+\nu_{\gamma})}{\Gamma(1-\nu_{\gamma})}.
\label{eqq11}
\end{equation}
We begin by considering the case $2\nu_{\gamma}>\tau$. Then
$$
{\tilde f}(\lambda)=2\lambda^{2\nu_{\gamma}}{\widetilde{\tilde f}}(\lambda),
$$
where
$$
\matrix{{\widetilde{\tilde f}}(\lambda)=
&=&
\cos A^*\cos A^{**}-\Frac{a_1}{l\lambda^{1}}\sin(A^*+A^{**})\hfill  \cr
&&+\Frac{1}{2\lambda^{2}}\left(4\left(b_2-\Frac{a_2}{l^2}\right)
\cos A^*\cos A^{**}+ 2\left(\Frac{a_1^2}{l^2}-b_1^2\right)\sin A^*\sin A^{**}\right)\hfill \cr
&&+\Frac{\nu_{\gamma}\alpha b_1}{\lambda^{3}}\sin(A^{**}-A^{*})+O\left(\Frac{1}{\lambda^{4}}\right)\hfill \cr
&&+C\left[\Frac{1}{2\lambda^{2\nu_{\gamma}-\tau}}\cos^2 A^{**}
+ \Frac{1}{2\lambda^{1+2\nu_{\gamma}-\tau}}
\left(-\Frac{a_1}{l}\sin(2A^{**})
+(\tau-1)\omega\cos^2 A^{**}\right)\right.\hfill \cr
&&+\Frac{1}{2\lambda^{2+2\nu_{\gamma}-\tau}}\Bigg(2\left(b_2-\Frac{a_2}{l^2}\right)
\cos^2 A^{**}+\left(\Frac{a_1^2}{l^2}-b_1^2\right)
\sin^2 A^{**}\hfill \cr
&&\left.-(\tau-1)\omega\Frac{a_1}{l}\sin(2A^{**})+\Frac{(\tau-1)(\tau-2)}{2}\omega^2\cos^2 A^{**}\Bigg)\right]\hfill \cr
&&+O(\lambda^{\tau-3-2\nu_{\gamma}}).\hfill \cr
}
$$
Then
\begin{equation}
{\widetilde{\tilde f}}(\lambda)=f_0(\lambda)+\Frac{f_1(\lambda)}{\lambda^{2\nu_{\gamma}-\tau}}
+\Frac{f_2(\lambda)}{\lambda}+\Frac{f_3(\lambda)}{\lambda^{1+2\nu_{\gamma}-\tau}}+\Frac{f_4(\lambda)}{\lambda^2}
+\Frac{f_5(\lambda)}{\lambda^{2+2\nu_{\gamma}-\tau}}+O\left(\frac{1}{\lambda^3}\right),
\label{027}
\end{equation}
where
\begin{equation}
f_0(\lambda)=\cos A^*\cos A^{**},
\label{032}
\end{equation}
\begin{equation}
f_1(\lambda)=\Frac{C}{2}\cos^2 A^{**},
\label{032s}
\end{equation}
\begin{equation}
f_2(\lambda)=-\Frac{a_1}{l}\sin(A^*+A^{**}),
\label{032d}
\end{equation}
\begin{equation}
f_3(\lambda)=\Frac{C}{2}\left(-\Frac{a_1}{l}\sin(2A^{**})+(\tau-1)\omega\cos^2 A^{**}\right),
\label{eqqq3}
\end{equation}
\begin{equation}
f_4(\lambda)=\Frac{1}{2}\left(4\left(b_2-\Frac{a_2}{l^2}\right)
\cos A^*\cos A^{**}+ 2\left(\Frac{a_1^2}{l^2}-b_1^2\right)\sin A^*\sin A^{**}\right),
\label{eqqq4}
\end{equation}
\begin{equation}
\matrix{f_5(\lambda)&=&\Frac{C}{2}\Bigg(2\left(b_2-\Frac{a_2}{l^2}\right)\cos^2 A^{**}+\left(\Frac{a_1^2}{l^2}-b_1^2\right)
\sin^2 A^{**}\hfill \cr
&&-(\tau-1)\omega\Frac{a_1}{l}\sin(2A^{**})+\Frac{(\tau-1)(\tau-2)}{2}\omega^2\cos^2 A^{**}\Bigg).\cr}
\label{eqqq5}
\end{equation}
Note that $f_0, f_1, f_2, f_3, f_4$ and $f_5$ remain bounded in the strip $-\alpha_0\leq
\Re(\lambda)\leq 0$.

\noindent
{\bf Step 2.}
We now investigate the zeros of $f_0$. It follows from (\ref{032}) that $f_0$ has two families of zeros, which we denote by $\lambda_{1, k}^0$
and $\lambda_{2, k}^0$, respectively.
$$
f_0(\lambda)=0 \Leftrightarrow \cos A^*\cos A^{**}=0.
$$
This leads to
$$
\left\{\matrix{iq\lambda+\nu_{\gamma}\frac{\pi}{2}-\frac{\pi}{4}=k\pi+\frac{\pi}{2} \hbox{ or }\hfill &\cr
iq\lambda-\nu_{\gamma}\frac{\pi}{2}-\frac{\pi}{4}=k\pi+\frac{\pi}{2}\hfill &\cr}\right.
$$
which immediately yields the following expressions for the roots:
$$
\matrix{\lambda_{1, k}^0&=&\frac{1}{iq}(k\pi-\nu_{\gamma}\frac{\pi}{2}+\frac{3}{4}\pi),\quad k\in\Z,\hfill \cr
\lambda_{2, k}^0&=&\frac{1}{iq}(k\pi+\nu_{\gamma}\frac{\pi}{2}+\frac{3}{4}\pi),\quad k\in\Z.\hfill \cr}
$$
By Rouch\'e's Theorem, we conclude that ${\tilde f}$ has infinitely many simple zeros in $\Sigma$, which can be divided into two families,
denoted by $\lambda_{1, k}$ and $\lambda_{2, k}$, respectively. More precisely, for $|k|\geq k_0$, where $k_0$ is sufficiently large, we have
\begin{equation}
\lambda_{1, k}=\lambda_{1, k}^0+\varepsilon_{1, k},
\label{030eaa}
\end{equation}
\begin{equation}
\lambda_{2, k}=\lambda_{2, k}^0+\varepsilon_{2, k}.
\label{030e}
\end{equation}
{\bf Step 3. Asymptotic behavior of $\varepsilon_{1, k}$.}
Using (\ref{030eaa}), we get
\begin{equation}
\matrix{\cos(qi\lambda_{1, k}-\nu_{\gamma}\frac{\pi}{2}-\frac{\pi}{4})&=&(-1)^{k+1}\sin (-\nu_{\gamma}\frac{\pi}{2}+iq\varepsilon_{1, k})\hfill  \cr
\cos(qi\lambda_{1, k}+\nu_{\gamma}\frac{\pi}{2}-\frac{\pi}{4})&=&(-1)^{k+1}\sin (iq\varepsilon_{1, k}).\hfill  \cr}
\label{eqq14}
\end{equation}
Substituting (\ref{eqq14}) into (\ref{032}), using that
${\tilde f}(\lambda_{1, k})=0$, we get
\begin{equation}
{\tilde f}(\lambda_{1, k})=
-iq\varepsilon_{1, k}\sin(\nu_{\gamma}\pi)+\Frac{C\sin(\nu_{\gamma}\pi)^2}{2(\frac{1}{iq}k\pi)^{2\nu_{\gamma}-\tau}}+O\left(\Frac{1}{k}\right)=0.
\label{eqq15}
\end{equation}
The previous equation has one solution
\begin{equation}
\varepsilon_{1, k}=-\Frac{\tilde C}{q^{\tau+1}}\Frac{\sin(\nu_{\gamma}\pi)}{(k\pi)^{2\nu_{\gamma}-\tau}}i^{1-\tau}+O\left(\Frac{1}{k}\right),
\label{eqq16}
\end{equation}
where
$$
{\tilde C}=\Frac{2^{2\nu_{\gamma}+1}\rho}{(1-\gamma)}\Frac{\Gamma(1+\nu_{\gamma})}{\Gamma(1-\nu_{\gamma})}.
$$
From (\ref{eqq16}) we have in that case
$|k|^{2\nu_{\gamma}-\tau}\Re\lambda_{1, k}\sim \upsilon$ with
$$
\upsilon=-\Frac{\tilde C}{q^{\tau+1}}\Frac{\sin(\nu_{\gamma}\pi)\sin\frac{\pi}{2}\tau}{\pi^{2\nu_{\gamma}-\tau}}.
$$
{\bf Step 4. Asymptotic behavior of $\varepsilon_{2, k}$.}
Using (\ref{030e}), we get
\begin{equation}
\matrix{\cos(qi\lambda_{2, k}-\nu_{\gamma}\frac{\pi}{2}-\frac{\pi}{4})&=&(-1)^{k+1}\sin (iq\varepsilon_{2, k})\hfill  \cr
&=&(-1)^{k+1}\left[iq\varepsilon_{2, k}+\frac{i q^3}{6}\varepsilon_{2, k}^3+O(\varepsilon_{2, k}^3)\right]\hfill  \cr
\cos(qi\lambda_{2, k}+\nu_{\gamma}\frac{\pi}{2}-\frac{\pi}{4})&=&(-1)^{k+1}\sin(\nu_{\gamma}\pi+iq\varepsilon_{2, k})\hfill  \cr
&=&(-1)^{k+1}\left[\sin(\nu_{\gamma}\pi)+iq\cos(\nu_{\gamma}\pi)\varepsilon_{2, k}
+\frac{q^2}{2}\sin(\nu_{\gamma}\pi)\varepsilon_{2, k}^2
\right.\hfill  \cr
&&\left.+\frac{iq^3}{6}\cos(\nu_{\gamma}\pi)\varepsilon_{2, k}^3+O(\varepsilon_{2, k}^3)\right]\hfill  \cr
\sin(A^*+A^{**})&=&-\sin(\nu_{\gamma}\pi+2iq\varepsilon_{2, k})\hfill  \cr
&=&-\left[\sin(\nu_{\gamma}\pi)+2 i q\cos(\nu_{\gamma}\pi)\varepsilon_{2, k}
+2 q^2\sin(\nu_{\gamma}\pi)\varepsilon_{2, k}^2
\right.\hfill  \cr
&&\left.+\frac{4 i q^3}{3}\cos(\nu_{\gamma}\pi)\varepsilon_{2, k}^3+O(\varepsilon_{2, k}^3)\right].\hfill  \cr
}
\label{ddff}
\end{equation}
Substituting (\ref{ddff}) into (\ref{032}), using that
${\tilde f}(\lambda_{2, k})=0$, we get
\begin{equation}
{\tilde f}(\lambda_{2, k})=
iq\varepsilon_{2, k}\sin(\nu_{\gamma}\pi)+\Frac{a_1\sin(\nu_{\gamma}\pi)}{k\pi}+O(\varepsilon_{2, k}^2)+O\left(\Frac{1}{k^2}\right)+O\left(\Frac{\varepsilon_{2, k}}{k}\right)=0.
\label{hh1}
\end{equation}
The previous equation has one solution
\begin{equation}
\varepsilon_{2, k}=i\frac{a_1}{q k\pi}+O(\varepsilon_{2, k}^2)+O\left(\Frac{1}{k^2}\right)+O\left(\Frac{\varepsilon_{2, k}}{k}\right).
\label{hh2}
\end{equation}

\noindent
We can write
\begin{equation}
\lambda_{2, k}=\frac{1}{iq}(k\pi+\nu_{\gamma}\frac{1}{2}+\frac{3}{4}\pi)
+i\frac{a_1}{q k\pi}+{\tilde\varepsilon}_{2, k}.
\label{hh3}
\end{equation}
where ${\tilde\varepsilon}_{2, k}=o(1/k)$.
Substituting (\ref{hh3}) into (\ref{027}), we get
\begin{equation}
\matrix{{\tilde f}(\lambda_{2, k})=
iq{\tilde\varepsilon}_{2, k}\sin(\nu_{\gamma}\pi)
-\Frac{a_1(\frac{\nu_{\gamma}}{2}+\frac{3}{4})}{k^2\pi}\sin(\nu_{\gamma}\pi)
+\left(2\Frac{a_1^2}{l^2}-b_1^2\right)\Frac{\cos(\nu_{\gamma}\pi)}{(\frac{1}{iq}k\pi)^2}
+O({\tilde\varepsilon}_{2, k}^2)+O\left(\Frac{1}{k^3}\right)\hfill \cr
+O\left(\Frac{{\tilde\varepsilon}_{2, k}}{k}\right)=0.\cr}
\label{mml1}
\end{equation}
The previous equation gives
\begin{equation}
{\tilde\varepsilon}_{2, k}=-\Frac{a_1(\frac{\nu_{\gamma}}{2}+\frac{3}{4})}{iq k^2\pi}
+\left(2\Frac{a_1^2}{l^2}-b_1^2\right)\Frac{\cos(\nu_{\gamma}\pi)}{\sin(\nu_{\gamma}\pi)}\Frac{iq}{(k\pi)^2}
+O({\tilde\varepsilon}_{2, k}^2)+O\left(\Frac{1}{k^3}\right)+O\left(\Frac{{\tilde\varepsilon}_{2, k}}{k}\right).
\label{hh5}
\end{equation}
We can write
\begin{equation}
\lambda_{2, k}=\frac{1}{iq}(k\pi+\nu_{\gamma}\frac{1}{2}+\frac{3}{4}\pi)
+i\frac{a_1}{q k\pi}-\Frac{a_1(\frac{\nu_{\gamma}}{2}+\frac{3}{4})}{iq k^2\pi}
+\left(2\Frac{a_1^2}{l^2}-b_1^2\right)\Frac{\cos(\nu_{\gamma}\pi)}{\sin(\nu_{\gamma}\pi)}\Frac{iq}{(k\pi)^2}
+{\widetilde{\tilde\varepsilon}}_{2, k}.
\label{hh3hh3}
\end{equation}
where ${\widetilde{\tilde\varepsilon}}_{2, k}=o(1/k^2)$.
Substituting (\ref{hh3hh3}) into (\ref{027}), we get
\begin{equation}
\matrix{{\tilde f}(\lambda_{2, k})=
iq{\widetilde{\tilde\varepsilon}}_{2, k}\sin(\nu_{\gamma}\pi)+C\left[\Frac{a_1^2}{2(\frac{1}{iq}k\pi)^{2\nu_{\gamma}-\tau}(k\pi)^2}-\Frac{a_1^2}{l(\frac{1}{iq}k\pi)^{1+2\nu_{\gamma}-\tau}(k\pi)}\right.\hfill &\cr
\left.+\Frac{\frac{a_1^2}{l^2}-b_1^2}{2(\frac{1}{iq}k\pi)^{2+2\nu_{\gamma}-\tau}}\right]
+O\left(\Frac{1}{k^{3}}\right)=0, &\cr}
\label{hh5hh5}
\end{equation}
where $C$ is defined in (\ref{eqq11}). then
\begin{equation}
{\widetilde{\tilde\varepsilon}}_{2, k}=-\Frac{{\tilde C}}{8}\Frac{\alpha^2 q^{3-\tau}}{\sin(\nu_{\gamma}\pi)}\Frac{i^{1-\tau}}{(k\pi)^{2+2\nu_{\gamma}-\tau}}+O\left(\Frac{1}{k^{3}}\right).
\label{eqq12}
\end{equation}
From (\ref{hh5}) we have in that case
$|k|^{2+2\nu_{\gamma}-\tau}\Re\lambda_{2, k}\sim \upsilon$ with
$$
\upsilon=-\Frac{{\tilde C}}{8}q^{3-\tau}\Frac{1}{\pi^{2+2\nu_{\gamma}-\tau}}\Frac{\sin\frac{\pi}{2}\tau}{\sin(\nu_{\gamma}\pi)}.
$$

\noindent
{\bf Case $2\nu_{\gamma}<\tau$}.
Then
$$
{\tilde f}(\lambda)=C\lambda^{\tau}{\widetilde{\tilde f}}(\lambda),
$$
where
$$
\matrix{{\widetilde{\tilde f}}(\lambda)
&=&\cos^2 A^{**}
+ \Frac{1}{\lambda}
\left(-\Frac{a_1}{l}\sin(2A^{**})
+(\tau-1)\omega\cos^2 A^{**}\right)\hfill \cr
&&+\Frac{1}{\lambda^2}\Bigg(2\left(b_2-\Frac{a_2}{l^2}\right)
\cos^2 A^{**}+\left(\Frac{a_1^2}{l^2}-b_1^2\right)
\sin^2 A^{**}\hfill \cr
&&-(\tau-1)\omega\Frac{a_1}{l}\sin(2A^{**})+\Frac{(\tau-1)(\tau-2)}{2}\omega^2\cos^2 A^{**}\Bigg)\hfill \cr
&&\Frac{1}{C}\left[\Frac{2}{\lambda^{\tau-2\nu_{\gamma}}}\cos A^{*}\cos A^{**}
-\Frac{2a_1}{l\lambda^{1+\tau-2\nu_{\gamma}}}
\sin(A^*+A^{**})\right.\hfill  \cr
&&+\left.\Frac{1}{\lambda^{2+\tau-2\nu_{\gamma}}}\left(4\left(b_2-\Frac{a_2}{l^2}\right)
\cos A^*\cos A^{**}+ 2\left(\Frac{a_1^2}{l^2}-b_1^2\right)\sin A^*\sin A^{**}\right)\right]\hfill \cr
&&+O\left(\Frac{1}{\lambda^3}\right). \cr}
$$
Then
\begin{equation}
{\tilde f}(\lambda)=f_0(\lambda)+\Frac{f_1(\lambda)}{\lambda^{\tau-2\nu_{\gamma}}}+\Frac{f_2(\lambda)}{\lambda}
+\Frac{f_3(\lambda)}{\lambda^{1+\tau-2\nu_{\gamma}}}+\Frac{f_4(\lambda)}{\lambda^{2}}+\Frac{f_5(\lambda)}{\lambda^{2+\tau-2\nu_{\gamma}}}+O\left(\frac{1}{\lambda^3}\right),
\label{027h}
\end{equation}
where
\begin{equation}
f_0(\lambda)=\cos^2 A^{**},
\label{032h}
\end{equation}
\begin{equation}
f_1(\lambda)=\Frac{2}{C}\cos A^{*}\cos A^{**},
\label{032dh}
\end{equation}
\begin{equation}
f_2(\lambda)=-\Frac{a_1}{l}\sin(2A^{**})+(\tau-1)\omega\cos^2 A^{**},
\label{032dh1}
\end{equation}
\begin{equation}
f_3(\lambda)=-\Frac{2a_1}{C l}\sin(A^*+A^{**}),
\label{032sh}
\end{equation}
\begin{equation}
\matrix{f_4(\lambda)&=&\Bigg(2\left(b_2-\Frac{a_2}{l^2}\right)
\cos^2 A^{**}+\left(\Frac{a_1^2}{l^2}-b_1^2\right)
\sin^2 A^{**}\hfill \cr
&&-(\tau-1)\omega\Frac{a_1}{l}\sin(2A^{**})+\Frac{(\tau-1)(\tau-2)}{2}\omega^2\cos^2 A^{**}\Bigg),\hfill \cr}
\label{032dh2}
\end{equation}
\begin{equation}
f_5(\lambda)=\Frac{1}{C}\left(4\left(b_2-\Frac{a_2}{l^2}\right)
\cos A^*\cos A^{**}+ 2\left(\Frac{a_1^2}{l^2}-b_1^2\right)\sin A^*\sin A^{**}\right).
\label{032dh3}
\end{equation}
Note that $f_0,f_1, f_2, f_3, f_4$ and $f_5$ remain bounded in the strip $-\alpha_0\leq
\Re(\lambda)\leq 0$.

\noindent
{\bf Step 2.} We look at the roots of $f_0$. From
(\ref{032h}), $f_0$ has two families of roots that we denote
$\lambda_{1, k}^0$ and $\lambda_{2, k}^0$.
$$
f_0(\lambda)=0 \Leftrightarrow \cos^2 A^{**}=0.
$$
Then
$$
\matrix{\lambda_{1, k}^0&=&
\frac{1}{iq}(k\pi+\nu_{\gamma}\frac{\pi}{2}+\frac{3}{4}\pi),\quad k\in\Z,\hfill\cr
\lambda_{2, k}^0&=&\frac{1}{iq}(k\pi+\nu_{\gamma}\frac{\pi}{2}+\frac{3}{4}\pi),\quad k\in\Z.\hfill\cr}
$$
By Rouch\'e's Theorem, we conclude that ${\tilde f}$ has infinitely many simple zeros in $\Sigma$, which can be divided into two families,
denoted by $\lambda_{1, k}$ and $\lambda_{2, k}$, respectively. More precisely, for $|k|\geq k_0$, where $k_0$ is sufficiently large, we have
\begin{equation}
\lambda_{1, k}=\frac{1}{iq}(k\pi+\nu_{\gamma}\frac{\pi}{2}+\frac{3}{4}\pi)+\varepsilon_{1, k},
\label{eqq45}
\end{equation}
\begin{equation}
\lambda_{2, k}=\frac{1}{iq}(k\pi+\nu_{\gamma}\frac{\pi}{2}+\frac{3}{4}\pi)+\varepsilon_{2, k}.
\label{eqq46}
\end{equation}
We have
\begin{equation}
\matrix{\cos(l\lambda_{j, k}-\nu_{\gamma}\frac{\pi}{2}-\frac{\pi}{4})&=&(-1)^{k+1}\sin (l\varepsilon_{j, k})=(-1)^{k+1}(l\varepsilon_{j, k}+O(\varepsilon_{j, k}^3)), \hfill \cr
\cos(l\lambda_{j, k}+\nu_{\gamma}\frac{\pi}{2}-\frac{\pi}{4})&=&(-1)^{k+1}\sin (\nu_{\gamma}\pi+l\varepsilon_{j, k}) \hfill \cr
&=&(-1)^{k+1}(\sin(\nu_{\gamma}\pi)+l\varepsilon_{j, k}\cos(\nu_{\gamma}\pi)+O(\varepsilon_{j, k}^2)), \hfill \cr}
\label{eqq47}
\end{equation}
\begin{equation}
\matrix{\sin(l\lambda_{j, k}-\nu_{\gamma}\frac{\pi}{2}-\frac{\pi}{4})&=&(-1)^{k}\cos (l\varepsilon_{j, k})=(-1)^{k}(1+O(\varepsilon_{j, k}^2)), \hfill \cr
\sin(l\lambda_{j, k}+\nu_{\gamma}\frac{\pi}{2}-\frac{\pi}{4})&=&(-1)^{k}\cos (\nu_{\gamma}\pi+l\varepsilon_{j, k})\hfill \cr
&=&(-1)^{k}(\cos (\nu_{\gamma}\pi)  -l\varepsilon_{j, k}\sin(\nu_{\gamma}\pi)+O(\varepsilon_{j, k}^2))\hfill \cr}
\label{eqq48}
\end{equation}
\begin{equation}
\Frac{1}{(\lambda_{j, k})^{r}}=\Frac{l^{r}}{(k\pi)^r}\left(1-r\Frac{\nu_{\gamma}\frac{\pi}{2}+\frac{3}{4}}{k}-r\Frac{l\varepsilon_{j, k}}{k\pi}+O\left(\Frac{\varepsilon_{j, k}}{k^2}\right)+O\left(\Frac{1}{k^2}\right)\right).
\label{eqq49}
\end{equation}
Substituting (\ref{eqq47}), (\ref{eqq48}) and (\ref{eqq49}) into (\ref{027h}), using that ${\widetilde{\tilde f}}(\lambda)=0$, we have
$$
\matrix{l^2\varepsilon_{j, k}^2+\frac{2a_1 l}{k\pi}\left(1-\frac{\frac{\nu_{\gamma}}{2}+\frac{3}{4}}{k}\right)\varepsilon_{j, k}
+\frac{2 l}{C}\sin(\nu_{\gamma}\pi)\frac{l^{\tau-2\nu_{\gamma}}}{(k\pi)^{(\tau-2\nu_{\gamma})}}
\left(1-\frac{(\tau-2\nu_{\gamma})(\frac{\nu_{\gamma}}{2}+\frac{3}{4})}{k}\right)\varepsilon_{j, k}\hfill & \cr
+\left(\frac{a_1^2}{l^2}-b_1^2\right)\frac{l^2}{(k\pi)^2}
+\frac{2a_1}{l C}\sin(\nu_{\gamma}\pi)\frac{l^{1+\tau-2\nu_{\gamma}}}{(k\pi)^{1+\tau-2\nu_{\gamma}}}+O(\varepsilon_{j, k}^2)+O\left(\Frac{1}{k^2}\right)=0.\hfill & \cr}
$$
The previous equation has two solutions
\begin{equation}
\matrix{\varepsilon_{1, k}&=&-\frac{2l^{\tau-2\nu_\gamma-1}\sin\nu_\gamma\pi}{C(k\pi)^{\tau-2\nu_\gamma}}+o\bigl(k^{-(\tau-2\nu_{\gamma})}\bigr),\hfill  \cr
\varepsilon_{2, k}&=&-\frac{a_1}{l\pi k}+\frac{{\tilde C} b_1^2}{2l^{\tau-1}\sin\nu_\gamma\,(k\pi)^{2+2\nu_{\gamma}-\tau}}+o\bigl(k^{-(2+2\nu_{\gamma}-\tau)}\bigr)\hfill  \cr
&=&-\frac{a_1}{l\pi k}+\frac{{\tilde C} \alpha^2}{2l^{\tau-3}\sin\nu_\gamma\pi\,(k\pi)^{2+2\nu_{\gamma}-\tau}}+o\bigl(k^{-(2+2\nu_{\gamma}-\tau)}\bigr).\hfill  \cr}
\label{eqq51}
\end{equation}
From (\ref{eqq51}) we have in that case
$$
|k|^{\tau-2\nu_{\gamma}}\Re\lambda_{1, k}\sim \upsilon_1\ \hbox{ and }\ |k|^{2+2\nu_{\gamma}-\tau}\Re\lambda_{2, k}\sim \upsilon_2\ \hbox{ with }
$$
$$
\upsilon_1=-2\Frac{q^{\tau-1}}{\tilde C}\Frac{\sin(\nu_{\gamma}\pi)\sin\frac{\pi}{2}\tau}{\pi^{\tau-2\nu_{\gamma}}},\quad
\upsilon_2=-\Frac{{\tilde C}q^{3-\tau} \alpha^2}{2}\Frac{\sin\frac{\pi}{2}\tau}{\sin(\nu_{\gamma}\pi)}   \Frac{1}{\pi^{2+2\nu_{\gamma}-\tau}}
$$
{\bf Case $\tau=2\nu_{\gamma}$}
Then
$$
{\tilde f}(\lambda)=\lambda^{\tau}{\widetilde{\tilde f}}(\lambda),
$$
where
\begin{equation}
{\widetilde{\tilde f}}(\lambda)=
f_0(\lambda)+\frac{1}{\lambda}f_1(\lambda)+\frac{1}{\lambda^2}f_2(\lambda)+O\left(\Frac{1}{\lambda^3}\right).
\label{eqq90}
\end{equation}
where
\begin{equation}
f_0(\lambda)=2\cos(A^{**}+d)\cos A^{**}+C\cos^2 A^{**},\hspace{5cm}
\label{eqq92}
\end{equation}
\begin{equation}
f_1(\lambda)=-\frac{2a_1}{l}\sin(2 A^{**}+d)
+C\left[-\frac{a_1}{l}\sin(2 A^{**})+(\tau-1)\omega\cos^2 A^{**}\right],\hspace{0.3cm}
\label{eqq93}
\end{equation}
\begin{equation}
\matrix{\hspace{0.8cm} f_2(\lambda)=& 4\left(b_2-\frac{a_2}{l^2}\right)\cos(A^{**}+d)\cos A^{**}
+2\left(\frac{a_1^2}{l^2}-b_1^2\right)\sin(A^{**}+d)\sin A^{**}\hfill  \cr
&+C\Bigg[2\left(b_2-\frac{a_2}{l^2}\right)\cos^2 A^{**}
+\left(\frac{a_1^2}{l^2}-b_1^2\right)\sin^2 A^{**}\hfill  &\cr
&\qquad-(\tau-1)\omega\frac{a_1}{l}\sin(2 A^{**})
+\frac{(\tau-1)(\tau-2)}{2}\omega^2\cos^2 A^{**}\Bigg].\hfill  \cr}
\label{eqq94}
\end{equation}
where $d=\frac{\tau}{2}\pi$.
Note that $f_0,f_1$ and $f_2$ remain bounded in the strip $-\alpha_0\leq\Re(\lambda)\leq 0$.

\noindent
{\bf Step 2.} We look at the roots of $f_0$. From
(\ref{eqq92}), $f_0$ has two families of roots that we denote
$\lambda_{1, k}^0$ and $\lambda_{2, k}^0$.
$$
\matrix{f_0(\lambda)=0 &\Leftrightarrow & 2\cos(A^{**}+d)\cos A^{**}+c\cos^2 A^{**}=0\hfill \cr
&\Leftrightarrow & \cos A^{**}=0 \hbox{ or } 2\cos(A^{**}+d)+c\cos A^{**}=0\hfill \cr}
$$
Firstly, we have
$$
 \cos A^{**}=0\Leftrightarrow \lambda_{1, k}^0=
\frac{1}{iq}(k\pi+\nu_{\gamma}\frac{\pi}{2}+\frac{3}{4}\pi),\quad k\in\Z.
$$
Using Rouch\'e Theorem, we obtain
$$
\matrix{\lambda_{1, k}&=&
\frac{1}{iq}(k\pi+\nu_{\gamma}\frac{\pi}{2}+\frac{3}{4}\pi)+\frac{i a_1}{q\pi}\frac{1}{k}
+\frac{1}{k^2}\Bigg[-\frac{i a_1}{q\pi^2}\left(\frac{3\pi}{4}+\frac{\pi\tau}{4}\right)\hfill  \cr
&& \qquad+\frac{i q b_1^2}{2\pi^2\sin(\pi\tau/2)}
\left(2\cos\frac{\pi\tau}{2}+\frac{{\tilde C}}{q^\tau}\cos\frac{\pi\tau}{2}-i\frac{{\tilde C}}{q^\tau}\sin\frac{\pi\tau}{2}\right)
\Bigg]+O(k^{-3})\hfill \cr
&=&\left[-\frac{1}{q}(k\pi+\nu_{\gamma}\frac{\pi}{2}+\frac{3}{4}\pi)+\frac{ a_1}{q\pi}\frac{1}{k}+\left(-\frac{a_1}{q\pi^2}\left(\frac{3\pi}{4}+\frac{\pi\tau}{4}\right)
+\frac{q b_1^2}{2\pi^2\sin(\pi\tau/2)}
\left(2\cos\frac{\pi\tau}{2}+\frac{{\tilde C}}{q^\tau}\cos\frac{\pi\tau}{2}\right)\right)\frac{1}{k^2}\right]i\hfill  \cr
&&+\Frac{q^{1-\tau} b_1^2 {\tilde C}}{2\pi^2}\frac{1}{k^2}+O(k^{-3})\hfill  \cr}
$$
Then
$$
|k|^{2}\Re\lambda_{1, k}\sim \upsilon_1=-\Frac{q^{3-\tau}\alpha^2 {\tilde C}}{8\pi^2}.
$$
Secondly, we have
$$
2\cos(A^{**}+d)\cos A^{**}+c\cos^2 A^{**}=0\Leftrightarrow e^{-2q\lambda}=\Frac{-(2+r)}{(2-r)\sin d-i(2+r)\cos d},
$$
where $r=\frac{{\tilde C}}{q^{\tau}}$.

Let us set $\lambda=x+iy$. Then, we have
$$
\left\{\matrix{e^{-2qx}=\Frac{2+r}{\sqrt{4+r^2+4r\cos \tau}}, \hfill \cr
-\frac{4}{2-\gamma}y=2k\pi-\Frac{\pi}{2}+\theta,\quad k\in \Z, \hfill \cr}\right.
$$
where $\theta$ is such that
$$
\left\{\matrix{\cos\theta&=&\Frac{(2+r)\cos\nu_\gamma\pi}{\sqrt{4+r^2+4r\cos2\nu_\gamma\pi}},\hfill\cr
\sin\theta&=&\Frac{(2-r)\sin\nu_\gamma\pi}{\sqrt{4+r^2+2r\cos2\nu_\gamma\pi}}.\hfill\cr}\right.
$$
Hence
$$
\left\{\matrix{x=-\frac{2-\gamma}{4}\ln\Frac{2+r}{\sqrt{4+r^2+4r\cos 2\nu_\gamma\pi}}, \hfill \cr
y=-\frac{2-\gamma}{4}(2k\pi-\Frac{\pi}{2}+\theta),\quad k\in \Z. \hfill \cr}\right.
$$
Now with the help of Rouch\'e's Theorem, we conclude that
$$
\lambda_{2, k}^0=-i\frac{2-\gamma}{4}(2k\pi-\Frac{\pi}{2}+\theta)-\frac{2-\gamma}{4}\ln\Frac{2+r}{\sqrt{4+r^2+4r\cos 2\nu_\gamma\pi}}+O(k^{-1}),\quad k\in\Z.
$$
Now, setting ${\tilde U}_k=(\lambda_{j, k}^0-{\cal P}) U_k$, where $U_k$ is
a normalized eigenfunction associated to $\lambda_{j, k}$. We then have
$$
\matrix{\|(\lambda_{2, k}^0-{\cal P})^{-1}\|_{{\cal L}({\cal H})}=\Sup_{U\in {\cal H}, U\not=0} \Frac{\|(\lambda_{2, k}^0-{\cal P})^{-1}U\|_{{\cal H}}}{\|U\|_{{\cal H}}}
&\geq &\Frac{\|(\lambda_{2, k}^0-{\cal P})^{-1}{\tilde U}_k\|_{{\cal H}}}{\|{\tilde U}_k\|_{{\cal H}}}\hfill \cr
&\geq &\Frac{\|U_k\|_{{\cal H}}}{\|(\lambda_{2, k}^0-{\cal P}) U_k\|_{{\cal H}}}.\hfill \cr}
$$
Hence, by Lemma \ref{ll12}, we deduce that
$$
\|(\lambda_{2, k}^0-{\cal P})^{-1}\|_{{\cal L}({\cal H})}\geq c \left\{\matrix{|k|^{2+2\nu_{\gamma}-\tau}\hfill & \hbox{ if }2\nu_{\gamma}>\tau,\hfill \cr
|k|^{2+2\nu_{\gamma}-\tau}\hfill & \hbox{ if }2\nu_{\gamma}<\tau.\hfill \cr
|k|^{2}\hfill & \hbox{ if }2\nu_{\gamma}=\tau.\hfill \cr}\right.
$$
Consequently, the semigroup $e^{t {\cal P}}$ is not exponentially stable, which completes the proof.

\hfill$\diamondsuit$\\

\section{Polynomial stability (for $\omega\not=0$)}
To establish polynomial energy decay, we invoke the following theorem.
\begin{theorem}[\cite{boto}]
Assume that ${\cal P}$ is the generator of a strongly continuous semigroup of contractions $(e^{t{\cal P}})_{t\geq 0}$
on a Hilbert space ${\cal X}$. If $i\R\subset\rho({\cal P})$. Then for a fixed $\wp> 0$ the following conditions are equivalent
\begin{itemize}
\item[1)] $\Sup_{\beta\in \R}\|(is I-{\cal P})^{-1}\|_{{\cal L}({\cal X})}=O(|s|^{\wp})$.
\item[2)] $\|e^{t{\cal P}}U_0\|_{\cal X}\leq \Frac{C}{t^{\frac{1}{\wp}}} \|U_0\|_{D({\cal P})}\quad\forall t>0,\ U_0\in D({\cal P})$,
for some $C> 0$.
\end{itemize}
\label{thm2}
\end{theorem}
Our main result is the following.
\begin{theorem}
The semigroup ${S_{\cal P}(t)}_{t\geq 0}$ associated with system $(P')$ is polynomially stable, i.e., there exists
a constant $C> 0$ such that
$$
{\cal E}(t)=\|S_{\cal P}(t)U_0\|_{\cal H}^2\leq\left\{\matrix{
\Frac{c}{ t^{  \frac{2}{2+2\nu_\gamma-\tau }}}\|U_0\|_{D({\cal P})}^2 &\hbox{ if }2\nu_\gamma>\tau,\hfill \cr
\Frac{c}{ t^{\frac{2}{2+2\nu_\gamma-\tau}}}\|U_0\|_{D({\cal P})}^2 &\hbox{ if }2\nu_\gamma<\tau,\hfill \cr
\Frac{c}{ t}\|U_0\|_{D({\cal P})}^2 &\hbox{ if }2\nu_\gamma=\tau.\hfill \cr}\right.
$$
Moreover, the rate of energy decay is optimal for general initial data in $D({\cal P})$.

\label{thm225}
\end{theorem}
{\bf Proof}\\
In section 3, we have proved that the first condition in Theorem \ref{thm2} is satisfied. Now,we need to show that
\begin{equation}
\Sup_{|\lambda|\geq 1}\Frac{1}{\lambda^{\wp}}\|(i\lambda I-{\cal P})^{-1}\|_{{\cal H}}< \infty,
\label{v35}
\end{equation}
for some positive constant $\wp$.

Let $\lambda\in \R$ and $F=(f_1, f_2 ,f_3, f_4, f_5)^{T}\in {\cal H}$. Consider the resolvent equation:
\begin{equation}
(i\lambda I-{\cal P})U=F.
\label{eqq110}
\end{equation}
Equivalently, we have
\begin{equation}
\left\{\matrix{i\lambda u-{\tilde u}=f_1,  \hfill & \cr
i\lambda {\tilde u}-(x^{\gamma}u_x)_x+\alpha v=f_2, \hfill & \cr
i\lambda v-{\tilde v}=f_3\hfill & \cr
i\lambda {\tilde u}-(x^{\gamma}v_x)_x+\alpha u=f_4, \hfill & \cr
(i\lambda+\varsigma^2+\omega)\varphi-{\tilde u}(0)\vartheta(\varsigma)=f_5,}\right.
\label{eqq56}
\end{equation}
with boundary conditions
\begin{equation}
\left\{\matrix{u(1)=v(0)=v(1)=0, \hfill & \cr
 (x^{\gamma}u_x)(0)=\zeta\Int_{-\infty}^{\infty}\vartheta(\varsigma)\varphi(\varsigma)\, d\varsigma,\hfill &\cr}\right.
\label{eqq57}
\end{equation}
where $F=(f_1, f_2, f_3, f_4, f_5)^{T}\in {\cal H}$.
From $(\ref{eqq56})_1, (\ref{eqq56})_2$ and $(\ref{eqq56})_3, (\ref{eqq56})_4$, we have
\begin{equation}
\left\{\matrix{-\lambda^2u- (x^{\gamma}u_x)_x+\alpha v=(f_2+i\lambda f_1), \hfill & \cr
-\lambda^2v- (x^{\gamma}v_x)_x+\alpha u=(f_4+i\lambda f_3). \hfill & \cr}\right.
\label{eqq58}
\end{equation}
Let us set
\begin{equation}
\left\{\matrix{ \phi=u+v, \hfill \cr
 {\psi}=u-v.
\hfill &\cr}\right.
\label{eqq59}
\end{equation}
Then, we obtain
 \begin{equation}
 \left\{\matrix{(\lambda^2-\alpha)\phi-(x^{\gamma}\phi_x)_x=-(i\lambda f_1+f_2+i\lambda f_3+f_4),\hfill &\cr
(\lambda^2+\alpha)\psi-(x^{\gamma}\psi_x)_x=-(i\lambda f_1+f_2-i\lambda f_3+f_4).
\hfill &\cr}\right.
\label{eqq60}
 \end{equation}
The solution of equations (\ref{eqq60}) is given by
\begin{equation}
\left\{\matrix{\phi(x)=A\Phi_{+}(x)+B\Phi_{-}(x)\hfill &\cr
\quad +\Frac{\pi}{(2-\gamma)\sin\nu_{\gamma}\pi}\Int_{0}^{x}(i\lambda f_1+f_2+i\lambda f_3+f_4)(\Phi_{+}(s)\Phi_{-}(x)-\Phi_{+}(x)\Phi_{-}(s))\, ds,&\cr
\psi(x)={\tilde A}\Phi_{++}+{\tilde B}\Phi_{--}\hfill &\cr
\quad +\Frac{\pi}{(2-\gamma)\sin\nu_{\gamma}\pi}\Int_{0}^{x}(i\lambda f_1+f_2-i\lambda f_3-f_4)(\Phi_{++}(s)\Phi_{--}(x)-\Phi_{++}(x)\Phi_{--}(s))\, ds, &\cr
}\right.
\label{eqq61}
\end{equation}
where $\Phi_{+}, \Phi_{-}, \Phi_{++}$ and $\Phi_{--}$ are defined by
\begin{equation}
\left\{\matrix{\Phi_{+}(x):=x^{\frac{1-\gamma}{2}}J_{\nu_\gamma}\left(\frac{2}{2-\gamma}{\tilde\lambda}
x^{\frac{2-\gamma}{2}}\right),\
\Phi_{-}(x):=x^{\frac{1-\gamma}{2}}J_{-\nu_\gamma}\left(\frac{2}{2-\gamma}{\tilde\lambda}
x^{\frac{2-\gamma}{2}}\right),\hfill &\cr
\Phi_{++}(x):=x^{\frac{1-\gamma}{2}}J_{\nu_\gamma}\left(\frac{2}{2-\gamma}{\tilde{\tilde\lambda}}
x^{\frac{2-\gamma}{2}}\right),\
\Phi_{--}(x):=x^{\frac{1-\gamma}{2}}J_{-\nu_\gamma}\left(\frac{2}{2-\gamma}{\tilde{\tilde\lambda}}
x^{\frac{2-\gamma}{2}}\right),\hfill &\cr}\right.
\label{eqq101}
\end{equation}
where
$$
\nu_\gamma=\Frac{1-\gamma}{2-\gamma},\ {\tilde\lambda}=\sqrt{\lambda^2-\alpha},\ {\tilde{\tilde\lambda}}=\sqrt{\lambda^2+\alpha}.
$$
Then
$$
\left\{\matrix{u(x)=\Frac{1}{2}(\phi+\psi),\hfill &\cr
v(x)=\Frac{1}{2}(\phi-\psi).\hfill &\cr}\right.
$$
Then
\begin{equation}
\matrix{u'(x)=\Frac{1}{2}(A\Phi'_{+}(x)+B\Phi'_{-}(x)+{\tilde A}\Phi'_{++}+{\tilde B}\Phi'_{--})\hfill &\cr
+\Frac{\pi}{2(2-\gamma)\sin\nu_{\gamma}\pi}\left[
\Int_{0}^{x}(i\lambda f_1+f_2+i\lambda f_3+f_4)(\Phi_{+}(s)\Phi'_{-}(x)-\Phi'_{+}(x)\Phi_{-}(s))\, ds\right.\hfill &\cr
\left.+\Int_{0}^{x}(i\lambda f_1+f_2-i\lambda f_3-f_4)(\Phi_{++}(s)\Phi'_{--}(x)-\Phi'_{++}(x)\Phi_{--}(s))\, ds\right]\hfill &\cr}
\label{eqq64}
\end{equation}
The substitution of $\varphi$ given by $(\ref{eqq56})_3$ into $(\ref{eqq56})$ gives us
\begin{equation}
(x^{\gamma}u_x)(0)=\rho (i\lambda+\omega)^{\tau-1} {\tilde u}(0)
+\zeta\Int_{-\infty}^{+\infty}\Frac{\vartheta(\varsigma)f_5(\varsigma)}{i\lambda+\varsigma^{2}+\omega}\, d\varsigma.
\label{eqq62}
\end{equation}
Moreover, from $(\ref{eqq60})_1$, we have
$$
{\tilde u}(0)=i\lambda u(0)-f_1(0).
$$
Then, the condition (\ref{eqq57}) become
\begin{equation}
(x^{\gamma}u_x)(0)-\rho i\lambda(i\lambda+\omega)^{\tau-1} u(0)
=-\rho(i\lambda+\omega)^{\tau-1}f_1(0)+\zeta\Int_{-\infty}^{+\infty}\Frac{\vartheta(\varsigma)f_5(\varsigma)}{i\lambda+\varsigma^{2}+\omega}\, d\varsigma.
\label{eqq63}
\end{equation}
From (\ref{eqq63}), (\ref{eqq64}) and (\ref{eqq57}), we conclude that
\begin{equation}
\matrix{(1-\gamma){\tilde c}_{\nu_\gamma, 0}^{+}A+(1-\gamma){\widetilde{\tilde c}}_{\nu_\gamma, 0}^{+}{\tilde A}-\rho i\lambda(i\lambda+\omega)^{\tau-1}\left[{\tilde c}_{\nu_\gamma, 0}^{-}B
+{\widetilde{\tilde c}}_{\nu_\gamma, 0}^{-}{\tilde B}\right]\hfill & \cr
=-\rho(i\lambda+\omega)^{\tau-1}f_1(0)+\zeta\Int_{-\infty}^{+\infty}\Frac{\vartheta(\varsigma)f_3(\varsigma)}{i\lambda+\varsigma^{2}+\omega}\, d\varsigma,\hfill & \cr}
\label{eqq67}
\end{equation}
\begin{equation}
B{\tilde c}_{\nu_\gamma, 0}^{-}-{\tilde B}{\widetilde{\tilde c}}_{\nu_\gamma, 0}^{-}=0,
\label{eqq68}
\end{equation}
\begin{equation}
\matrix{A\Phi_{+}(1)+B\Phi_{-}(1)+{\tilde A}\Phi_{++}(1)+{\tilde B}\Phi_{--}(1)\hfill & \cr
+\Frac{\pi}{(2-\gamma)\sin\nu_{\gamma}\pi}\left[\Int_{0}^{1}(i\lambda f_1+f_2+i\lambda f_3+f_4)(\Phi_{+}(s)\Phi_{-}(1)-\Phi_{+}(1)\Phi_{-}(s))\, ds\right. &\cr
\left.+\Int_{0}^{1}(i\lambda f_1+f_2-i\lambda f_3-f_4)(\Phi_{++}(s)\Phi_{--}(1)-\Phi_{++}(1)\Phi_{--}(s))\, ds\right]=0, & \cr}
\label{eqq69}
\end{equation}
\begin{equation}
\matrix{A\Phi_{+}(1)+B\Phi_{-}(1)-{\tilde A}\Phi_{++}(1)-{\tilde B}\Phi_{--}(1)\hfill & \cr
+\Frac{\pi}{(2-\gamma)\sin\nu_{\gamma}\pi}\left[\Int_{0}^{1}(i\lambda f_1+f_2+i\lambda f_3+f_4)(\Phi_{+}(s)\Phi_{-}(1)-\Phi_{+}(1)\Phi_{-}(s))\, ds\right. &\cr
\left.-\Int_{0}^{1}(i\lambda f_1+f_2-i\lambda f_3-f_4)(\Phi_{++}(s)\Phi_{--}(1)-\Phi_{++}(1)\Phi_{--}(s))\, ds\right]=0, & \cr}
\label{eqq70}
\end{equation}
where
$$
\matrix{{\tilde c}_{\nu_\gamma, m}^{+}=c_{\nu_\gamma, m}^{+}\left(\Frac{2}{2-\gamma}{\tilde\lambda}\right)^{2m+\nu_\gamma},\quad
{\tilde c}_{\nu_\gamma, m}^{-}=c_{\nu_\gamma, m}^{-}\left(\Frac{2}{2-\gamma}{\tilde\lambda}\right)^{2m-\nu_\gamma}\hfill & \cr
{\widetilde{\tilde c}}_{\nu_\gamma, m}^{+}=c_{\nu_\gamma, m}^{+}\left(\Frac{2}{2-\gamma}{\widetilde{\tilde\lambda}}\right)^{2m+\nu_\gamma},\quad
{\tilde c}_{\nu_\gamma, m}^{-}=c_{\nu_\gamma, m}^{-}\left(\Frac{2}{2-\gamma}{\widetilde{\tilde\lambda}}\right)^{2m-\nu_\gamma}\hfill & \cr}
$$
and
$$
\matrix{\Phi_{+}(1)=J_{\nu_{\gamma}}\left(\frac{2}{2-\gamma}{\tilde\lambda}\right),\qquad
\Phi_{-}(1)=J_{-\nu_{\gamma}}\left(\frac{2}{2-\gamma}{\tilde\lambda}\right),\hfill & \cr
\Phi_{++}(1)=J_{\nu_{\gamma}}\left(\frac{2}{2-\gamma}{\widetilde{\tilde\lambda}}\right),\qquad
\Phi_{--}(1)=J_{-\nu_{\gamma}}\left(\frac{2}{2-\gamma}{\widetilde{\tilde\lambda}}\right).\hfill & \cr}
$$
Using (\ref{eqq67}), (\ref{eqq68}), (\ref{eqq69}) and (\ref{eqq70}), a linear system in $A, {\tilde A}, B$ and ${\tilde B}$ is obtained
\begin{equation}
\pmatrix{r_{11}& r_{12}& r_{13}& r_{14}\cr
r_{21}&
r_{22}& r_{23}&r_{24}\cr
r_{31}& r_{32}& r_{33}&r_{34}\cr
r_{41}& r_{42}& r_{33}&r_{44}\cr}\pmatrix{A\cr {\tilde A}\cr B\cr {\tilde B}\cr}=\pmatrix{C_1\cr C_2\cr C_3\cr C_4\cr},
\label{eq111}
\end{equation}
where
$$
\matrix{r_{11}=(1-\gamma){\tilde c}_{\nu_\gamma, 0}^{+},\quad
r_{12}=(1-\gamma){\widetilde{\tilde c}}_{\nu_\gamma, 0}^{+},\hfill &\cr
r_{13}=-\rho i\lambda(i\lambda+\omega)^{\tau-1}{\tilde c}_{\nu_\gamma, 0}^{-},\quad
r_{14}=-\rho i\lambda(i\lambda+\omega)^{\tau-1}{\widetilde{\tilde c}}_{\nu_\gamma, 0}^{-},\hfill &\cr
r_{21}=0,\quad r_{22}=0,\quad r_{23}= {\tilde c}_{\nu_\gamma, 0}^{-},\quad
r_{24}=-{\widetilde{\tilde c}}_{\nu_\gamma, 0}^{-},\hfill &\cr
r_{31}= J_{\nu_\gamma}\left(\frac{2}{2-\gamma}{\tilde\lambda}\right),\
r_{32}= J_{\nu_\gamma}\left(\frac{2}{2-\gamma}{\widetilde{\tilde\lambda}}\right),\
r_{33}= J_{-\nu_\gamma}\left(\frac{2}{2-\gamma}{\tilde\lambda}\right),\
r_{34}= J_{-\nu_\gamma}\left(\frac{2}{2-\gamma}{\widetilde{\tilde\lambda}}\right),\hfill &\cr
r_{41}= J_{\nu_\gamma}\left(\frac{2}{2-\gamma}{\tilde\lambda}\right),\
r_{42}= -J_{\nu_\gamma}\left(\frac{2}{2-\gamma}{\widetilde{\tilde\lambda}}\right),\
r_{43}= J_{-\nu_\gamma}\left(\frac{2}{2-\gamma}{\tilde\lambda}\right),\
r_{44}= -J_{-\nu_\gamma}\left(\frac{2}{2-\gamma}{\widetilde{\tilde\lambda}}\right),\hfill &\cr}
$$
$$
\matrix{C_1&=&-\rho(i\lambda+\omega)^{\tau-1}f_1(0)+\zeta\Int_{-\infty}^{+\infty}\Frac{\vartheta(\varsigma)f_3(\varsigma)}{i\lambda+\varsigma^{2}+\omega}\, d\varsigma,\hfill \cr
C_2&=&0,\hfill \cr
C_3&=&-\Frac{\pi}{(2-\gamma)\sin\nu_{\gamma}\pi}\left[\Int_{0}^{1}(i\lambda f_1+f_2+i\lambda f_3+f_4)(\Phi_{+}(s)\Phi_{-}(1)-\Phi_{+}(1)\Phi_{-}(s))\, ds\right. \hfill\cr
&&\left.+\Int_{0}^{1}(i\lambda f_1+f_2-i\lambda f_3-f_4)(\Phi_{++}(s)\Phi_{--}(1)-\Phi_{++}(1)\Phi_{--}(s))\, ds\right],\hfill \cr
C_4&=&-\Frac{\pi}{(2-\gamma)\sin\nu_{\gamma}\pi}\left[\Int_{0}^{1}(i\lambda f_1+f_2+i\lambda f_3+f_4)(\Phi_{+}(s)\Phi_{-}(1)-\Phi_{+}(1)\Phi_{-}(s))\, ds\right. \hfill\cr
&&\left.-\Int_{0}^{1}(i\lambda f_1+f_2-i\lambda f_3-f_4)(\Phi_{++}(s)\Phi_{--}(1)-\Phi_{++}(1)\Phi_{--}(s))\, ds\right].\hfill\cr}
$$
Let the determinant of the linear system given in (\ref{eq111}) be denoted by $D$. Then
$$
\matrix{D=
-2\left[(1-\gamma){\tilde c}_{\nu_\gamma, 0}^{+} {\widetilde{\tilde c}}_{\nu_\gamma, 0}^{-}
J_{\nu_\gamma}\left(\frac{2}{2-\gamma}{\tilde\lambda}\right)
J_{-\nu_\gamma}\left(\frac{2}{2-\gamma}{\widetilde{\tilde\lambda}}\right)\right.\hfill &\cr
+(1-\gamma) {\tilde c}_{\nu_\gamma, 0}^{-} {\widetilde{\tilde c}}_{\nu_\gamma, 0}^{+}
J_{\nu_\gamma}\left(\frac{2}{2-\gamma}{\widetilde{\tilde\lambda}}\right)
J_{-\nu_\gamma}\left(\frac{2}{2-\gamma}{\tilde\lambda}\right)\hfill &\cr
\left.+2\rho i\lambda(i\lambda+\omega)^{\tau-1}{\tilde c}_{\nu_\gamma, 0}^{-}
{\widetilde{\tilde c}}_{\nu_\gamma, 0}^{-}
J_{\nu_\gamma}\left(\frac{2}{2-\gamma}{\tilde\lambda}\right)
J_{\nu_\gamma}\left(\frac{2}{2-\gamma}{\widetilde{\tilde\lambda}}\right)\right]\hfill &\cr
=
-\Frac{2(1-\gamma) c_{\nu_\gamma, 0}^{+} c_{\nu_\gamma, 0}^{-}}{\lambda}\left[\cos A^*\cos A^{**}-\Frac{a_1}{l\lambda}\sin(A^*+A^{**})\right.\hfill  &\cr
+\Frac{1}{2\lambda^{2}}\left(4\left(b_2-\Frac{a_2}{l^2}\right)
\cos A^*\cos A^{**}+ 2\left(\Frac{a_1^2}{l^2}-b_1^2\right)\sin A^*\sin A^{**}\right)\hfill &\cr
-\Frac{\nu_{\gamma}\alpha b_1}{\lambda^{3}}\sin(A^{**}-A^{*})+O\left(\Frac{1}{\lambda^{4}}\right)\hfill &\cr
+\Frac{2^{2\nu_{\gamma}+1}\rho}{(1-\gamma) q^{2\nu_{\gamma}}}\Frac{\Gamma(1+\nu_{\gamma})}{\Gamma(1-\nu_{\gamma})}\left[\Frac{i^{\tau}}{2\lambda^{2\nu_{\gamma}-\tau}}\cos^2 A^{**}
+ \Frac{1}{2\lambda^{1+2\nu_{\gamma}-\tau}}
\left(-\Frac{i^{\tau}a_1}{l}\sin(2A^{**})
+\Frac{(\tau-1)}{i^{1-\tau}}\omega\cos^2 A^{**}\right)\right.\hfill &\cr
+\Frac{1}{2\lambda^{2+2\nu_{\gamma}-\tau}}\Bigg(2i^{\tau}\left(b_2-\Frac{a_2}{l^2}\right)
\cos^2 A^{**}+i^{\tau}\left(\Frac{a_1^2}{l^2}-b_1^2\right)
\sin^2 A^{**}\hfill &\cr
\left.\left.-\Frac{(\tau-1)}{i^{1-\tau}}\omega\Frac{a_1}{l}\sin(2A^{**})+\Frac{(\tau-1)(\tau-2)}{2 i^{2-\tau}}\omega^2\cos^2 A^{**}\Bigg)\right]\right]\hfill &\cr
+O(\lambda^{\tau-3-2\nu_{\gamma}}).\hfill &\cr}
$$
where
$$
b_1=\Frac{q\alpha}{2},\ b_2=-\Frac{q^2\alpha^2}{8},\ q=\frac{2}{2-\gamma},\ l=q.
$$
As $D\not=0$ for all $\lambda\not=0$, then $A$ and $B$ are uniquely determined by (\ref{eq111}).

Now, we have that the following important Lemma.
\begin{lemma}
\begin{equation}
|D|\geq c |\lambda|^{-3+\tau-2\nu_{\gamma}} \hbox{ for large } \lambda.
\label{qqcc1}
\end{equation}
\end{lemma}
{\bf Proof of Lemma \ref{qqcc1}.} See Appendix.

\noindent
Now
$$
B=\frac{r_{11}r_{24}r_{42}S-r_{12}r_{24}r_{31}T-r_{31}r_{32}C_1}{D}, \qquad
{\tilde B}=\frac{r_{11}r_{23}r_{32}S+r_{12}r_{23}r_{31}T-r_{31}r_{32}C_1}{D}.
$$
$$
A=\frac{r_{12}r_{34}r_{23}S+2r_{13}r_{24}r_{32}S+r_{24}r_{12}r_{33}T+r_{24}r_{33}r_{32}C_1}{D},
$$
$$
{\tilde A}=\frac{-r_{11}r_{33}r_{24}T+2r_{13}r_{24}r_{31}T-r_{23}r_{11}r_{34}S- r_{23}r_{34}r_{31}C_1}{D},
$$
where
$$
S=\Frac{C_3+C_4}{2},\quad T=\Frac{C_3-C_4}{2}.
$$
Considering only the dominant terms of $\lambda$, the following is obtained:
$$
\matrix{|D||B|, |D||{\tilde B}|&\leq & c_1|\lambda|^{-\frac{3}{2}}+c_2|\lambda|^{\tau-2}\leq \left\{\matrix{c_3|\lambda|^{-\frac{3}{2}} \hbox{ if }\tau\leq \frac{1}{2}\hfill &\cr
c_3|\lambda|^{\tau-2}\hbox{ if }\tau\geq \frac{1}{2}\hfill &\cr}\right.\hfill \cr
|D||A|, |D||{\tilde A}|&\leq & c_1|\lambda|^{-\frac{3}{2}}+c_2|\lambda|^{-\frac{3}{2}+\tau-2\nu_{\gamma}}+c_3|\lambda|^{-2+\tau-\nu_{\gamma}}\leq
{\tilde c}_3|\lambda|^{-\frac{3}{2}},\hfill \cr}
$$
Then, we conclude that
$$
 |A|,\ |{\tilde A}|\leq c |\lambda|^{\frac{3}{2}-\tau+2\nu_{\gamma}},
$$
$$
 |B|,\ |{\tilde B}|\leq \left\{\matrix{c |\lambda|^{\frac{3}{2}-\tau+2\nu_{\gamma}}\hbox{ if }\hbox{ if }\tau\leq \frac{1}{2},\hfill \cr
c |\lambda|^{1+2\nu_{\gamma}}\hbox{ if }\hbox{ if }\tau\geq \frac{1}{2}.\hfill \cr}\right.
$$
Then
$$
\|u\|_{L^{2}(0, 1)},\ \|v\|_{L^{2}(0, 1)}\leq c |\lambda|^{1+2\nu_{\gamma}-\tau}\left(\|f_1\|_{H_{0, \gamma}^{1}(0, 1)}+\|f_2\|_{L^{2}(0, 1)}+\|f_3\|_{L^{2}(-\infty, +\infty)}\right).
$$
Using $(\ref{eqq56})_1$ and $(\ref{eqq56})_3$, we get
$$
\|{\tilde u}\|_{L^{2}(0, 1)},\ \|{\tilde v}\|_{L^{2}(0, 1)}\leq c |\lambda|^{2+2\nu_{\gamma}-\tau}\left(\|f_1\|_{H_{0, \gamma}^{1}(0, 1)}+\|f_2\|_{L^{2}(0, 1)}+\|f_3\|_{L^{2}(-\infty, +\infty)}\right).
$$
From (\ref{eqq101}) and (\ref{e7112}), we have
$$
\left\{\matrix{x^{\gamma/2}\Phi'_{+}(x)&=&(\frac{1-\gamma}{2}+\frac{(2-\gamma)\nu_\gamma}{2})x^{-1/2}J_{\nu_\gamma}\left(\frac{2}{2-\gamma}{\tilde\lambda} x^{\frac{1-\gamma}{2}}\right)
-{\tilde\lambda} x^{\frac{1-\gamma}{2}}J_{1+\nu_\gamma}\left(\frac{2}{2-\gamma}{\tilde\lambda} x^{\frac{1-\gamma}{2}}\right),\hfill \cr
x^{\gamma/2}\Phi'_{-}(x)&=&(\frac{1-\gamma}{2}-\frac{(2-\gamma)\nu_\gamma}{2})x^{-1/2}J_{-\nu_\gamma}\left(\frac{2}{2-\gamma}{\tilde\lambda} x^{\frac{1-\gamma}{2}}\right)
-{\tilde\lambda} x^{\frac{1-\gamma}{2}}J_{1-\nu_\gamma}\left(\frac{2}{2-\gamma}{\tilde\lambda} x^{\frac{1-\gamma}{2}}\right).\hfill \cr}\right.
$$
Then from (\ref{eqq64}), we can get
$$
\|x^{\gamma/2}u_x\|_{L^{2}(0, 1)},\ \|x^{\gamma/2}v_x\|_{L^{2}(0, 1)}\leq c |\lambda|^{2+2\nu_{\gamma}-\tau}\left(\|f_1\|_{H_{0, \gamma}^{1}(0, 1)}+\|f_2\|_{L^{2}(0, 1)}+\|f_3\|_{L^{2}(-\infty, +\infty)}\right),
$$
Moreover
\begin{equation}
\zeta\Int_{-\infty}^{+\infty}(\varsigma^{2}+\omega)|\varphi(\varsigma)|^{2}\, d\varsigma\leq \|U\|_{{\cal H}}\|F\|_{{\cal H}}.
\label{qqww}
\end{equation}
Since $\omega> 0$, we have
$$
\|\varphi\|_{L^{2}(-\infty, \infty)}^2\leq \Int_{-\infty}^{+\infty}(\varsigma^{2}+\omega)|\varphi(\varsigma)|^{2}\, d\varsigma\leq c\|U\|_{{\cal H}}\|F\|_{{\cal H}}.
$$
Thus, we conclude that
$$
\|(i\lambda I-{\cal P})^{-1}\|_{{\cal H}}\leq \left\{\matrix{c |\lambda|^{2+2\nu_{\gamma}-\tau} & \hbox{ as } |\lambda|\rightarrow \infty& \hbox{ if }\tau\not=2\nu_{\gamma},  \hfill &\cr
c |\lambda|^{2} \ &\hbox{ as } |\lambda|\rightarrow \infty  &\hbox{ if }\tau=2\nu_{\gamma}. \hfill \cr}\right.
$$
The conclusion then follows by applying Theorem \ref{thm2}.
\section{Polynomial stability (for $\omega=0$)}
To prove polynomial decay, we use the following theorem.
\begin{theorem}[\cite{batty}]
Let $S(t)$ be a bounded $C_0$-semigroup on a Hilbert space ${\cal X}$ with generator ${\cal P}$.
Assume that $\sigma({\cal P})\cap i\R=\{0\}$ and that there exist
$\vartheta\geq 1$ and $\upsilon > 0$ such that
$$
\|(is I-{\cal P})^{-1}\|_{{\cal L}({\cal X})}=\left\{\matrix{O(|s|^{-\vartheta}),\quad & s\rightarrow 0, \hfill & \cr
O(|s|^{\upsilon}),\quad & |s|\rightarrow \infty. \hfill & \cr}\right.
$$
Then there exist constants $C, t_0> 0$ such that for all $t\geq t_0$ and $U_0\in D({\cal P})\cap Ran({\cal P})$
we have
$$
\|e^{{\cal P}t}U_0\|^{2}\leq C\Frac{1}{t^{\frac{2}{\varsigma}}}
\|U_0\|_{D({\cal A)}\cap Ran({\cal P})}^{2}\quad (\hbox{ with }\varsigma=\max\{\vartheta, \upsilon\}).
$$
\label{thm2kkk}
\end{theorem}
Our main results are the following.
\begin{theorem}
The semigroup ${S_{{\cal P}}(t)}_{t\geq 0}$ is polynomially stable and
$$
{\cal E}(t)=\|S_{{\cal P}_1}(t)U_0\|_{\cal H}^2\leq \Frac{1}{t^{2/(2+2\nu_{\gamma}-\tau)}}  \|U_0\|_{D({\cal P}_1)\cap Ran({\cal P})}^2.
$$
\label{thm225kkk}
\end{theorem}
{\bf Proof of Theorem \ref{thm225kkk}}
We will study the behavior of the resolvent operator $(i\lambda-{\cal P})^{-1}$ as $|\lambda|\rightarrow 0$.
For $\lambda$ near zero, we have
$$
|D|\geq c> 0.
$$
Moreover
$$
|A|, |{\tilde A}|, |B|, |{\tilde B}|\leq c |\lambda|^{\frac{\tau-2}{2}}(\|f_1\|_{L^2(0, 1)}+\|f_2\|_{L^2(0, 1)}+\|f_3\|_{L^2(0, 1)}+\|f_4\|_{L^2(0, 1)}+\|f_5\|_{L^2(-\infty, +\infty)})
$$
Then
$$
\|u\|_{L^{2}(0, 1)}, \|v\|_{L^{2}(0, 1)}\leq c |\lambda|^{\frac{\tau-2}{2}}(\|f_1\|_{L^2(0, 1)}+\|f_2\|_{L^2(0, 1)}+\|f_3\|_{L^2(0, 1)}+\|f_4\|_{L^2(0, 1)}+\|f_5\|_{L^2(-\infty, +\infty)}).
$$
$$
\|{\tilde u}\|_{L^{2}(0, 1)}, \|{\tilde v}\|_{L^{2}(0, 1)}\leq c (\|f_1\|_{L^2(0, 1)}+\|f_2\|_{L^2(0, 1)}+\|f_3\|_{L^2(0, 1)}+\|f_4\|_{L^2(0, 1)}+\|f_5\|_{L^2(-\infty, +\infty)}).
$$
Moreover from $(\ref{eqq56})_5$
$$
\varphi=\Frac{{\tilde u}(0)\vartheta(\varsigma)+f_5(\varsigma)}{i\lambda+\varsigma^2}.
$$
Then
$$
\matrix{\|\varphi\|_{L^2(\R)}^2&\leq & 2|{\tilde u}(0)|^2\Int_{\R}\Frac{\vartheta(\varsigma)^2}{|i\lambda+\varsigma^2|^2}\, d\varsigma+2\Int_{\R}\Frac{|f52(\varsigma)|^2}{|i\lambda+\varsigma^2|^2}\, d\varsigma\hfill \cr
&\leq & c |\lambda|^{\tau-2}|{\tilde u}(0)|^2+\Frac{1}{|\lambda|^2}\|f_5\|_{L^2(\R)}^{2}\hfill \cr}
$$
Hence
$$
\|U\|^2\leq C|\lambda|^{-2}\|F\|^2 \hbox{ as } \lambda\rightarrow 0.
$$
Finally,
$$
\|U\|\leq \left\{\matrix{c |\lambda|^{2+2\nu_{\gamma}-\tau}\|F\|_{\cal H}\hbox{ as } \lambda\rightarrow +\infty,\hfill & \cr
c|\lambda|^{-1}\|F\|_{\cal H} \hbox{ as } \lambda\rightarrow 0.\hfill & \cr}\right.
$$
The conclusion then follows by applying Theorem \ref{thm2kkk}.
\section*{Appendix}
{\bf Proof of Lemma}
Indeed, suppose (\ref{qqcc1}) was wrong. We consider three cases $\tau< 2 \nu_{\gamma}, \tau> 2 \nu_{\gamma}$ and $\tau=\nu_{\gamma}$. Then
$\exists \lambda_n$ such that $|\lambda_n|\rightarrow \infty$ with
\begin{equation}
|D| |\lambda_n|^{3-\tau+2\nu_{\gamma}}\rightarrow 0 \hbox{ as } n\rightarrow \infty.
\label{qqcc2}
\end{equation}
We can construct a proof by contradiction by simultaneously exploiting the
real and imaginary parts of $D(\lambda)$. This approach is particularly useful
here, since it takes into account the oscillatory behavior of $A^*$ and
$A^{**}$.

From the expression of $D$, namely $D=\frac{-2(1-\gamma)c_{\nu_{\gamma}, 0}^+ c_{\nu_{\gamma}, 0}^-}{\lambda}{\tilde D}$, it is sufficient to establish the lower bound
$$
{\tilde D}\geq K |\lambda|^{-(2+2\nu_{\gamma}-\tau)},
$$
for some constant $K$.

\noindent
We first consider the regime
$$
r=2\nu_\gamma-\tau>0.
$$
Our objective is to prove that there exist constants $K>0$ and
$\lambda_0>0$ such that
$$
|{\tilde D}(\lambda)|\geq K\lambda^{-(r+2)}
\qquad\hbox{ for all }\lambda\geq\lambda_0.
$$

\section*{The case $2\nu_\gamma>\tau$}

\subsection*{Setting up the contradiction}

Suppose, by contradiction, that the above lower bound is false. Then there
exists a sequence of real numbers $\lambda_n\to+\infty$ such that
$$
\lambda_n^{r+2}|{\tilde D}(\lambda_n)|\longrightarrow0.
$$
Since
$$
{\tilde D}=D_R+i\kappa E,\qquad \kappa\neq0,
$$
we necessarily have
\begin{equation}
\lambda_n^{r+2}D_R(\lambda_n)\longrightarrow0,
\label{1}
\end{equation}
and
\begin{equation}
\lambda_n^{r+2}E(\lambda_n)\longrightarrow0.
\label{2}
\end{equation}
In what follows, we omit the index $n$ for simplicity.

\subsection*{Consequence of the imaginary part}
Recall that
$$
\Frac{2}{\sin\frac{\tau}{2}\pi}\lambda^rE=
X^2+\frac{-2aXY+dX^2}{\lambda}
+\frac{pX^2+s_0Y^2-2daXY}{\lambda^2}+O\left(\Frac{1}{\lambda^{3}}\right),
$$
where
$$
X=\cos A^{**},\qquad Y=\sin A^{**},
\qquad s_0=a^2-b_1^2.
$$
$$
a=\frac{a_1}{l},\qquad d=(1-\tau)\Frac{\cos\frac{\tau}{2}\pi}{\sin\frac{\tau}{2}\pi}\omega.
$$
Assumption (\ref{2}) implies
$$
X^2+\frac{-2aXY+dX^2}{\lambda}
+\frac{pX^2+s_0Y^2-2daXY}{\lambda^2}=o(\lambda^{-2}).
$$
Multiplying by $\lambda^2$, we obtain
$$
\lambda^2X^2+-2a\lambda XY+d\lambda X^2
+pX^2+s_0Y^2-2daXY=o(1).
$$
Set
$$
T=\lambda X.
$$
We obtain
\begin{equation}
\left(1+\frac{d}{\lambda}+\frac{p}{\lambda^2}\right)T^2-2a\left(1+\frac{d}{\lambda^2}\right)YT+s_0 Y^2=o(1).
\label{eqq80}
\end{equation}
The coefficient of $T^2$ tends to $1$, whereas $|Y|\leq 1$. It follows that the sequences $(T_n)$ is bounded. Thus,
$$
X=\Frac{T_n}{\lambda_n}\rightarrow 0, |Y_n|\rightarrow 1.
$$
Passing to the limit in (\ref{eqq80}), we obtain
$$
T^2-2aYT+s_0 Y^2=o(1).
$$
since $s_0=a^2-b_1^2$, this becomes
\begin{equation}
(\lambda X-aY)^2-b_1^2Y^2=o(1).
\label{eqq81}
\end{equation}
Thus,
\begin{equation}
(\lambda X-aY)^2-b_1^2Y^2\longrightarrow0.
\label{3}
\end{equation}
Since $X^2+Y^2=1$, this relation implies
$$
\lambda X=O(1).
$$
In particular,
$$
X\longrightarrow0,\qquad |Y|\longrightarrow1.
$$
After extracting a subsequence, we may assume that
$$
Y\longrightarrow\varepsilon,\qquad \varepsilon\in\{-1,1\}.
$$

\subsection*{Consequence of the real part}
We have
$$
\matrix{D_R=ZX-\frac a\lambda(WX+ZY)
+\frac1{\lambda^2}(p_0ZX+s_0WY)\hfill &\cr
+k\Frac{\cos\frac{\tau}{2}\pi}{2\lambda^r}
\left(X^2+\frac{-2aXY+{\tilde d}X^2}{\lambda}
+\frac{pX^2+s_0Y^2-2{\tilde d}aXY}{\lambda^2}+O\left(\Frac{1}{\lambda^{3}}\right)\right),\hfill &\cr
=ZX-\frac a\lambda(WX+ZY)
+\frac1{\lambda^2}(p_0ZX+s_0WY)\hfill &\cr
+k\Frac{\cos\frac{\tau}{2}\pi}{2\lambda^r}
\left(X^2-2a\frac{XY}{\lambda}
+\frac{s_0Y^2}{\lambda^2}+O\left(\Frac{1}{\lambda^{3}}\right)\right),\hfill &\cr}
$$
where
$$
{\tilde d}=(\tau-1)\Frac{\sin\frac{\tau}{2}\pi}{\cos\frac{\tau}{2}\pi}\omega,
$$
$$
Z=X\cos\delta-Y\sin\delta,
$$
and
$$
W=Y\cos\delta+X\sin\delta,
\qquad \delta=\pi\nu_\gamma.
$$
Since $X\to0$ and $Y\to\varepsilon$,
$$
Z\longrightarrow-\varepsilon\sin\delta,
\qquad
W\longrightarrow\varepsilon\cos\delta.
$$
Using (\ref{1}) and multiplying the expression for $D_R$ by $\lambda$, we obtain
$$
\lambda D_R
=Z(\lambda X-aY)-aWX+O(\lambda^{-1}).
$$
Since $X=O(\lambda^{-1})$, it follows that
$$
\lambda D_R=Z(\lambda X-aY)+o(1).
$$
Assumption (\ref{1}) therefore implies
$$
Z(\lambda X-aY)\longrightarrow0.
$$
However,
$$
Z\longrightarrow-\varepsilon\sin(\pi\nu_\gamma)\neq0.
$$
Consequently,
\begin{equation}
\lambda X-aY\longrightarrow0.
\label{5}
\end{equation}

\subsection*{Contradiction}
Comparing (\ref{3}) and (\ref{5}), we obtain
$$
-b_1^2Y^2\longrightarrow0.
$$
But $Y^2\longrightarrow1$. It follows that $b_1^2=0$, which contradicts
the assumption $b_1(=q\alpha/2)\neq0$. This completes the contradiction argument.
\subsection*{Conclusion}
Under the assumptions
$$
a_1\neq0,\qquad b_1\neq0,\qquad 0<\tau<2\nu_\gamma,
$$
and with $\kappa\neq0$, we obtain
$$
\exists K>0,\quad \exists\lambda_0>0:
\qquad
|{\tilde D}(\lambda)|\geq
\frac{K}{\lambda^{2+2\nu_\gamma-\tau}},
\qquad \lambda\geq\lambda_0.
$$

\section*{The case $2\nu_\gamma<\tau$}
Let us now consider the case
$$
2\nu_\gamma<\tau.
$$
In this regime, the memory term in the imaginary part of $D(\lambda)$ is no
longer decaying: it contains the factor $\lambda^{\tau-2\nu_\gamma}$. The
proof by contradiction must therefore be adapted slightly.

Under the assumption $b_1\neq0$, we can still obtain a lower bound of the
same type, with the appropriate exponent.
\subsection*{Statement of the lower bound}
Set
$$
\sigma=\tau-2\nu_\gamma.
$$
Since $2\nu_\gamma<\tau<1$, we have $0<\sigma<1$.

Under the assumptions
$$
b_1\neq0,\qquad \kappa\neq0,
$$
and assuming that all coefficients are real, we shall prove that there exist
constants $K>0$ and $\lambda_0>0$ such that
$$
|{\tilde D}(\lambda)|\geq
\frac{K}{\lambda^{2-\sigma}}
=
\frac{K}{\lambda^{2+2\nu_\gamma-\tau}},
\qquad \lambda\geq\lambda_0.
$$

\subsection*{Proof by contradiction}
Recall the decomposition
$$
{\tilde D}(\lambda)=D_R(\lambda)+i\kappa E(\lambda),
$$
where
$$
E(\lambda)=\sin(\tau\frac{\pi}{2})\frac{\lambda^\sigma}{2}
\left[
X^2+\frac{-2aXY+dX^2}{\lambda}
+\frac{pX^2+s_0Y^2-2daXY}{\lambda^2}
+O\left(\frac{1}{\lambda^3}\right)\right],
$$
with
$$
a=\frac{a_1}{l},\qquad s_0=a^2-b_1^2,
\qquad
X=\cos A^{**},\qquad Y=\sin A^{**}.
$$
Moreover,
$$
\matrix{D_R=ZX-\frac{a}{\lambda}(WX+ZY)
+\frac1{\lambda^2}(p_0ZX+s_0WY)\hfill & \cr
+k\Frac{\cos\frac{\tau}{2}\pi}{2\lambda^{-\sigma}}
\left(X^2+\frac{-2aXY+{\tilde d}X^2}{\lambda}
+\frac{pX^2+s_0Y^2-2{\tilde d}aXY}{\lambda^2}+O\left(\Frac{1}{\lambda^{3}}\right)\right).\hfill &\cr}
$$
Suppose, by contradiction, that the announced lower bound is false. Then
there exists a sequence $\lambda_n\to+\infty$ such that
$$
\lambda_n^{2-\sigma}|D(\lambda_n)|\longrightarrow0.
$$
Since $\kappa\neq0$, we obtain
$$
\lambda_n^{2-\sigma}D_R(\lambda_n)\longrightarrow0,
\qquad
\lambda_n^{2-\sigma}E(\lambda_n)\longrightarrow0.
$$

\subsection*{Step 1. Exploiting the imaginary part}
The second convergence gives
$$
\lambda_n^2X_n^2
-2a\lambda_nX_nY_n+d\lambda_nX_n^2
+pX_n^2+s_0Y_n^2-2daX_nY_n
\longrightarrow0.
$$
Set
$$
T_n=\lambda_nX_n.
$$
We obtain
\begin{equation}
\left(1+\frac d{\lambda_n}+\frac p{\lambda_n^2}\right)T_n^2
-2a\left(1+\frac d{\lambda_n^2}\right)T_nY_n
+s_0Y_n^2\longrightarrow0.
\label{6}
\end{equation}
The coefficient of $T_n^2$ tends to $1$, whereas $|Y_n|\leq1$.
It follows that the sequence $(T_n)$ is bounded. Thus,
$$
X_n=\frac{T_n}{\lambda_n}\longrightarrow0,
\qquad |Y_n|\longrightarrow1.
$$
Passing to the limit in (\ref{6}), we obtain
$$
T_n^2-2aT_nY_n+s_0Y_n^2\longrightarrow0.
$$
Since $s_0=a^2-b_1^2$, this becomes
\begin{equation}
(T_n-aY_n)^2-b_1^2Y_n^2\longrightarrow0.
\label{eqq104}
\end{equation}

\subsection*{Step 2. Exploiting the real part}
Multiplying $D_R$ by $\lambda_n$, we have
$$
\matrix{\lambda_nD_R
=Z_nT_n-a(W_nX_n+Z_nY_n)
+\frac1{\lambda_n}(p_0Z_nX_n+s_0W_nY_n)\hfill &\cr
k\Frac{\cos\frac{\tau}{2}\pi}{2\lambda^{-(1+\sigma)}}
\left(X^2+\frac{-2aXY+{\tilde d}X^2}{\lambda}
+\frac{pX^2+s_0Y^2-2{\tilde d}aXY}{\lambda^2}+O\left(\Frac{1}{\lambda^{3}}\right)\right).\hfill &\cr}
$$
Since $T_n$ is bounded and $X_n=O(\lambda_n^{-1})$, it follows that
$$
\lambda_nD_R
=Z_n(T_n-aY_n)+O(\lambda_n^{-1+\sigma}).
\label{eqq105}
$$
Now, since $2-\sigma>1$,
$$
\lambda_nD_R
=\lambda_n^{-(1-\sigma)}
\big(\lambda_n^{2-\sigma}D_R\big)
\longrightarrow0.
$$
On the other hand, the relation
$$
Z_n=X_n\cos(\pi\nu_\gamma)
-Y_n\sin(\pi\nu_\gamma)
$$
and the preceding convergences imply
$$
|Z_n|\longrightarrow\sin(\pi\nu_\gamma)>0.
$$
Therefore, relation (\ref{eqq105}) yields
\begin{equation}
T_n-aY_n\longrightarrow0.
\label{eqq91}
\end{equation}

\subsection*{Step 3. Contradiction}
Combining (\ref{eqq104}) and (\ref{eqq91}), we obtain
$$
-b_1^2Y_n^2\longrightarrow 0.
$$
However, $Y_n^2\longrightarrow1$. It follows that $b_1=0$, which
contradicts our assumption $b_1\neq0$. The proof by contradiction is
complete.
\subsection{Conclusion}
We have therefore established, in the regime $2\nu_\gamma<\tau$,
$$
|{\tilde D}(\lambda)|\geq
K\lambda^{-(2+2\nu_\gamma-\tau)}
$$
for all sufficiently large $\lambda$.

\noindent
\textbf{Remark.}
Both regimes $2\nu_\gamma>\tau$ and $2\nu_\gamma<\tau$ lead to the same
lower-bound formula:
$$
|{\tilde D}(\lambda)|\geq
K\lambda^{-(2+2\nu_\gamma-\tau)}.
$$
However, the proof uses different estimates for the imaginary part,
depending on the sign of $2\nu_\gamma-\tau$.

The critical case $2\nu_\gamma=\tau$ can be treated separately using the
same argument, with the factor $\lambda^0=1$.

\end{sloppypar}
\end{document}